\documentclass[11pt,fancybox]{article}
\usepackage{amsmath,amssymb}
\usepackage{graphicx}
\newcommand{\dd}[1]{\langle #1\rangle}
\newcommand{\eproof}{\hspace*{\fill}{$\blacksquare$}}
\newtheorem{theorem}{Theorem}
\newtheorem{lemma}{Lemma}
\newtheorem{corollary}{Corollary}

\begin{document}
\title{
On a Generalization of the Bipartite Graph $D(k,q)$
\footnote{
The work of the first two authors was supported by the
Natural Science Foundation of China (No. 61379004).
The work of the third author was supported by the Singapore Ministry of Education under Research Grant MOE2016-T2-2-014 and by NTU under Tier 1 grant RG143/14.}
}
\author{
$\begin{array}{c} \text{Xiaoyan Cheng, Yuansheng Tang}\\
\text{\{xycheng,\ ystang\}@yzu.edu.cn}\\
\text{School of Mathematical Sciences}\\
\text{Yangzhou University,
P.R.China}\\
\text{and Huaxiong Wang}\\
\text{HXWang@ntu.edu.sg}\\
\text{School of Physical and Mathematical Sciences}\\
\text{Nanyang Technological University, Singapore}
\end{array}$}

\maketitle
\begin{abstract}
In this paper, we deal with a generalization $\Gamma(\Omega,q)$ of the bipartite
graphs $D(k,q)$ proposed by Lazebnik and Ustimenko, where $\Omega$ is a set of binary sequences
that are adopted to index the entries of the vertices.
A few sufficient
conditions on $\Omega$ for $\Gamma(\Omega,q)$ to admit a variety of automorphisms are proposed.
A sufficient condition for $\Gamma(\Omega,q)$ to be edge-transitive is
proposed further.
A lower bound of the number of the connected components of $\Gamma(\Omega,q)$
is given by showing some invariants for the components.
For $\Gamma(\Omega,q)$, paths and cycles which contain vertices of some specified form are investigated in details.
Some lower bounds for the girth
of $\Gamma(\Omega,q)$ are then shown.
In particular, one can give very simple conditions on
the index set $\Omega$ so as to assure the generalized graphs $\Gamma(\Omega,q)$ to be a family of graphs with large girth.


\end{abstract}

\section{Introduction}
All graphs we consider in this paper are assumed to be simple, i.e. undirected, without loops and multiple edges.
For a graph $G$, its vertex set and edge set are denoted by $V(G)$ and $E(G)$, respectively.
The order of $G$ is the number of the vertices in $V(G)$.
The size of $G$ is the number of the edges in $E(G)$.
The degree of a vertex of $G$ is the number of vertices adjacent to it.
A graph is said $r$-regular if the degrees of all its vertices are equal to $r$.
A sequence of vertices in $V(G)$ is called a path in $G$ if neighboring vertices are adjacent and neighbors of each vertex are different.
If $G$ is connected, the distance between two distinct vertices in $V(G)$ is the length of the shortest path connecting them
and the diameter of $G$ is the greatest distance between the vertices in $V(G)$.
A path is called a cycle further if its length is not smaller than 3 and
it is still a path when the beginning vertex is moved to the end for any number of rounds.
Clearly, a vertex may appear in a path/cycle several times.
If $G$ contains a cycle, then the girth of $G$, denoted by $g=g(G)$ is the length of the shortest cycles in $G$.
In literature, graphs with large girth and a high degree of symmetry
have been known to be hard to construct and have turned out to be
useful in different problems in extremal graph theory, finite geometry, coding theory, cryptography,
communication networks and quantum computations(\cite{K.P.P.P.F.2004},\cite{U.2007},\cite{U.2009},\cite{Y.T.2011}).

Let $q$ be a prime power and $\mathbb{F}_q$ the finite field of
$q$ elements. For $k\geq 2$, in \cite{L.U.1995} Lazebnik and
Ustimenko proposed a bipartite graph, denoted by $D(k,q)$,
which is $q$-regular, edge-transitive and of large girth.
The bipartite graph $D(k,q)$ can be equivalently described as the following~\cite{L.V.2004}:
The vertex sets $L(k)$ and $P(k)$ of $D(k,q)$ are two copies of $\mathbb{F}_q^{k}$
such that two vertices $(l_1,l_2,\ldots,l_{k})\in L(k)$ and
$(p_1,p_2,\ldots,p_{k})\in P(k)$ are adjacent in $D(k,q)$ if and only if
\begin{align}
    l_2+p_2&=p_1l_1,\label{-101}\\
    l_3+p_3&=p_1l_2,\label{-102}
\end{align}
and, for $4\leq i\leq k$,
\begin{align}
    l_i+p_i=\left\{
    \begin{array}{ll}
    -p_{i-2}l_1,&\text{if }i\equiv 0\text{ or }1(\bmod 4),\\
    p_1l_{i-2},&\text{if }i\equiv 2\text{ or }3(\bmod 4).
    \end{array}
    \right.\label{-103}
\end{align}
Clearly, if we define $l'_i=(-1)^{\lfloor
i/4\rfloor}l_i$ and $p'_i=(-1)^{\lfloor i/4\rfloor}p_i$
for $i\geq 1$, then (\ref{-101}), (\ref{-102}) and (\ref{-103})
can also be expressed as
\begin{align}
    l'_2+p'_2&=p'_1l'_1,\label{-104}\\
    l'_3+p'_3&=p'_1l'_2,\label{-105}
\end{align}
and
\begin{align}
    l'_i+p'_i&=\left\{
    \begin{array}{ll}
    l'_1p'_{i-2},&\text{if }i\equiv 0\text{ or }1(\bmod 4),\\
    p'_1l'_{i-2},&\text{if }i\equiv 2\text{ or }3(\bmod 4),
    \end{array}
    \right.\label{-106}
\end{align}
respectively, where $4\leq i\leq k$.
The construction of $D(k,q)$ was motivated by attempts to generalize
the notion of the "affine part" of a generalized polygon. In fact, $D(2,q)$ and $D(3,q)$ ($q$ odd) are
exactly the affine parts of a regular generalized 3-gon and 4-gon, respectively~\cite{L.U.1993,L.U.W.1999}.
In~\cite{L.U.1995}, some automorphisms of $D(k,q)$ were given and then the girth of
$D(k,q)$ was shown to be at least
$k+4$. Since the length of any cycle in bipartite graph must be even, this lower
bound for the girth of $D(k,q)$ is indeed $k+5$ for odd $k$.
In~\cite{F.L.S.U.W.1995}, it was proved that the girth of $D(k,q)$ is equal to $k+5$
if $k$ is odd and $(k+5)/2$ divides $q-1$.
It was also conjectured further in~\cite{F.L.S.U.W.1995} that the girth of $D(k,q)$ is equal to $k+5$ for odd $k$ and all $q\geq 4$.
This conjecture was proved in~\cite{CCT.2014} when $(k+5)/2$ is a power of the characteristic of $\mathbb{F}_q$
and in~\cite{CCT.2016} when $(k+5)/2$ is the product of a factor of $q-1$ and a power of the characteristic of $\mathbb{F}_q$, respectively.
For $k\geq 6$, it was shown in~\cite{L.U.W.1995} that the graph $D(k,q)$ is disconnected and
has at least $q^{t-1}$ components (any two being isomorphic), where $t=\lfloor(k+2)/4\rfloor$, by showing some invariants which are
fixed in each component.
Especially, this implied that the components of $D(k,q)$
provide the best-known asymptotic lower bound for the greatest number of edges in graphs of their order and girth.
The components of $D(k,q)$ are further characterized in~\cite{L.U.W.1996} for odd $q$ and in~\cite{L.V.2004}
for even $q$ with $q\geq 4$, respectively.
The properties of $D(k,q)$ were further investigated in~\cite{U.2007,U.2009} when the finite field $\mathbb{F}_q$
is replaced by a commutative ring.
Especially,
a new family of ordinary graphs of large girth was constructed in~\cite{U.2007} by using the natural polarity of $D(k,q)$.

For $n\geq 1$, Lazebnik and Viglione constructed in~\cite{L.V.2002} a bipartite graph $G_n(q)$
whose vertex sets are $L(n+1)$ and $P(n+1)$
such that two vertices $(l_1,l_2,\ldots,l_{n+1})\in L(n+1)$ and
$(p_1,p_2,\ldots,p_{n+1})\in P(n+1)$ are adjacent in $G_n(q)$ if and only if
\begin{align}
    l_i+p_i&=p_1l_{i-1}, i=2,3,\ldots,n+1.
\end{align}
For $n\geq 3$ and $q\geq 3$, or $n=2$ and $q$ odd,
the graph $G_n(q)$ is semi-symmetric~\cite{L.V.2002}, i.e. edge-transitive but not vertex-transitive.
It is also shown in~\cite{L.V.2002} that the graph $G_n(q)$ is connected when $1\leq n\leq q-1$ and
disconnected when $n\geq q$, in which case it has $q^{n-q+1}$
components, each isomorphic to $G_{q-1}(q)$. We note that $G_n(q)$ is indeed a generalization of the graph
defined by Wenger in~\cite{W.1991}. Thus, $G_n(q)$ is also called Wenger graph in~\cite{F.U.2007,V.2008}, where the diameter
of $G_n(q)$ was shown to be $2n+2$.

Let $\mathcal{R}$ be an arbitrary commutative ring.
Some bipartite graphs defined by systems of equations over $\mathcal{R}$ were investigated in~\cite{L.W.2001}.
For $i\geq 2$, let $f_i: \mathcal{R}^{2i-2}\rightarrow \mathcal{R}$ be any given function.
For $n\geq 2$, Lazebnik and Woldar construct in~\cite{L.W.2001} a bipartite graph $B\Gamma_n=B\Gamma(\mathcal{R}; f_2,\ldots,f_n)$
whose vertex sets $L_n$ and $P_n$ are two copies of $\mathcal{R}^n$ such that two vertices $(l_1,l_2,\ldots,l_{n})\in L_n$ and
$(p_1,p_2,\ldots,p_{n})\in P_n$ are adjacent in $B\Gamma_n$ if and only if
\begin{align}
    l_i+p_i&=f_i(p_1,l_1,p_2,l_2,\ldots,p_{i-1},l_{i-1}), i=2,3,\ldots,n.
\end{align}
Some general properties of $B\Gamma_n$ were exhibited in~\cite{L.W.2001}.
When $\mathcal{R}$ is a finite field and the functions $f_i$ are monomials of $p_1$ and $l_1$, the graph $B\Gamma_n$ is also called a monomial graph.
For positive integers $k,m,k',m'$, it is proved in~\cite{D.L.V.2005} that
the monomial graphs $B\Gamma_2(\mathbb{F}_q;p_1^kl_1^m)$ and $B\Gamma_2(\mathbb{F}_{q};p_1^{k'}l_1^{m'})$
are isomorphic if and only if $\{\gcd(k,q-1),\gcd(m,q-1)\}=\{\gcd(k',q-1),\gcd(m',q-1)\}$ as multi-sets.
If $q$ is odd, it was proved in~\cite{D.F.W.2007} that any monomial graph $B\Gamma_3$ of girth at least eight
is isomorphic to a graph $B\Gamma_3(\mathbb{F}_{q}; p_1l_1,p_1^kl_1^{2k})$ for some positive integer $k$ coprime to $q$.
In particular, the positive integer $k$ can be restricted to be 1 further if the odd prime power $q$ is not greater than $10^{10}$ or of form
$q=p^{2^a3^b}$ for odd prime $p$ and nonnegative integers $a,b$~\cite{D.F.W.2007}.
It was then conjectured in~\cite{D.F.W.2007} that, for any odd prime power $q$, every monomial graph $B\Gamma_3$ of girth at least eight
is isomorphic to $\Gamma_3=B\Gamma_3(\mathbb{F}_{q}; p_1l_1,p_1l_1^{2})$.

In this paper, we deal with only a specialized subclass of the graphs $B\Gamma_n$ that can also be seen as generalizations of $D(k,q)$ or $G_n(q)$.
Let $\Omega\subset\{0,1\}^*$ be a finite set of some binary
sequences. Let $\eta$ denote the null sequence. Throughout this paper, we
always assume that any binary sequence is of finite length and
$\Omega$ satisfies the following condition:
\begin{description}
    \item[C1]
    $\eta\in \Omega$ and for any sequence
$\alpha\in\Omega\setminus\{\eta\}$ the sequence obtained from
$\alpha$ by deleting its last bit is still in $\Omega$.
\end{description}
Let $L(\Omega)$ and $R(\Omega)$ be two copies of
$\mathbb{F}_q^{|\Omega|+1}$. We will denote the vectors in
$L(\Omega)$ and $R(\Omega)$ by $[l]$ and $\dd{r}$ respectively so
that we can distinguish the origin of vectors in the union set
$L(\Omega)\cup R(\Omega)$.
The entries of vectors in $L(\Omega)\cup R(\Omega)$ are indexed by the
elements in $\Omega\cup\{*\}$, where $*$ is a symbol not in
$\Omega$. For $[l]\in L(\Omega)$ and $\dd{r}\in R(\Omega)$, $l_*$ and $r_*$ are also called the colors of $[l]$ and $\dd{r}$, respectively.
Let $[0]_x$ denote the vector $[l]$ satisfying $l_*=x$ and $l_{\alpha}=0$ for $\alpha\in\Omega$,
and $\dd{0}_x$ denote the vector $\dd{r}$ satisfying $r_*=x$ and $r_{\alpha}=0$ for $\alpha\in\Omega$, respectively.
Thus, $[0]_x$ and $\dd{0}_x$ are the vectors which differ from the all-zero vector at most on the color entry.
Let $\Gamma(\Omega,q)$ be the bipartite graph with
$V(\Gamma(\Omega,q))=L(\Omega)\cup R(\Omega)$ and
$E(\Gamma(\Omega,q))\subset L(\Omega)\times R(\Omega)$ such that
$[l]\in L(\Omega)$ and $\dd{r}\in R(\Omega)$ are adjacent in
$\Gamma(\Omega,q)$, i.e. $([l],\dd{r})\in E(\Gamma(\Omega,q))$, if
and only if
\begin{align}
\label{000}
    l_{\eta}+r_{\eta}=l_*r_*,
\end{align}
and
\begin{align}\label{001}
    &l_{\alpha 0}+r_{\alpha 0}=r_*l_{\alpha},\text{ for }\alpha 0\in \Omega,\\
\label{002}
    &l_{\beta 1}+r_{\beta 1}=l_*r_{\beta},\text{ for }\beta 1\in \Omega.
\end{align}
If we define $*0=*1=\eta$, the equation (\ref{000}) can also be included in
either (\ref{001}) or (\ref{002}). Clearly, the bipartite graph
$\Gamma(\Omega,q)$ is $q$-regular, and the all-zero vectors
$[0]_0\in L(\Omega)$ and $\dd{0}_0\in R(\Omega)$ are adjacent in
$\Gamma(\Omega,q)$.
For $k\geq 2$, let $U_k$ denote the set consisting of the first
$k-1$ elements in the following set
\begin{align}
\label{0009}    U=\{\eta, 0,1,01,10,010,101,0101,1010,\ldots\}.
\end{align}
Then, $\Gamma(U_k,q)$ is equivalent to $D(k,q)$. For positive
integer $n$, let $W_n$ denote the set of all-0 binary sequences of
length less than $n$. Then, $\Gamma(W_n,q)$ is equivalent to the
Wenger graph $G_n(q)$.

This paper is arranged as follows. In Section~\ref{sec002}, we show some
automorphisms of $\Gamma(\Omega,q)$. A sufficient condition for
$\Gamma(\Omega,q)$ to be edge-transitive is proved by using these
automorphisms. In Section~\ref{sec000}, we will show some invariants which take fixed values in each components
of $\Gamma(\Omega,q)$. A lower bound for the number of components of $\Gamma(\Omega,q)$
is given further.
Behaviors of the components under some natural projections are also investigated in this section.
In Section~\ref{sec003}, for any path starting at a vertex of form $[0]_x$,
the vertices on the path are explicitly expressed by their colors.
Though a similar result for the paths starting at a vertex of form $\dd{0}_x$ can be written down directly,
it is omitted there for brevity.
In Section~\ref{sec004}, conditions for the existence of some cycles in $\Gamma(\Omega,q)$ are shown. In particular,
some lower bounds for the girth
of $\Gamma(\Omega,q)$ are shown.

\section{Automorphisms of $\Gamma(\Omega,q)$}\label{sec002}
In this section, we show some automorphisms of $\Gamma(\Omega,q)$.
For any binary sequence $\alpha$, let $|\alpha|$ and $w(\alpha)$
denote its length and the number of its nonzero bits,
respectively. For $x,y\in \mathbb{F}_q$, let $\lambda_{x,y}$ denote the map
over $L(\Omega)\cup R(\Omega)$ such that, for $[l]\in L(\Omega)$
and $\dd{r}\in R(\Omega)$,
\begin{align*}
    (\lambda_{x,y}([l]))_*=xl_*, (\lambda_{x,y}(\dd{r}))_*=yr_*,
\end{align*}
and
\begin{align*}
    (\lambda_{x,y}([l]))_{\alpha}&=x^{w(\alpha)+1}y^{|\alpha|-w(\alpha)+1}l_{\alpha},\\
    (\lambda_{x,y}(\dd{r}))_{\alpha}&=x^{w(\alpha)+1}y^{|\alpha|-w(\alpha)+1}r_{\alpha}.
\end{align*}
\begin{lemma}
\label{lem:0001} If $x,y\in \mathbb{F}_q$ are nonzero, then $\lambda_{x,y}$
is an automorphism of $\Gamma(\Omega,q)$.
\end{lemma}
{\bf Proof:} Since $x,y\in \mathbb{F}_q$ are nonzero, one can see easily,
from (\ref{000}) to (\ref{002}) and the definition of
$\lambda_{x,y}$, that $\lambda_{x,y}$ is a bijective map over
$L(\Omega)$ as well as over $R(\Omega)$ and,
$(\lambda_{x,y}([l]),\lambda_{x,y}(\dd{r}))\in
E(\Gamma(\Omega,q))$ if and only if $([l],\dd{r})\in
E(\Gamma(\Omega,q))$. Hence, $\lambda_{x,y}$ is an automorphism of
$\Gamma(\Omega,q)$.\eproof
\medskip

For $\alpha\in\Omega$, let $S_{\Omega}(\alpha)$ denote the set of sequences $\beta\in\{0,1\}^*$ with
\begin{align*}
    \{\alpha01\beta,\alpha10\beta\}\cap\Omega\neq\emptyset.
\end{align*}

\begin{lemma}
\label{lem:0002} If $\alpha$ is a sequence in $\Omega$ with $S_{\Omega}(\alpha)\subset\Omega$,
then, for any $x\in \mathbb{F}_q$, there is an automorphism
$\theta_{x,\alpha}$ of $\Gamma(\Omega,q)$ such that, for any
$[l]\in L(\Omega)$ and $\dd{r}\in R(\Omega)$,
\begin{align}
    (\theta_{x,\alpha}([l]))_{\alpha}=l_{\alpha}+x,\
    (\theta_{x,\alpha}(\dd{r}))_{\alpha}=r_{\alpha}-x \label{400},
\end{align}
and
\begin{align}
\label{401}
    (\theta_{x,\alpha}([l]))_{\gamma}=l_{\gamma},\
    (\theta_{x,\alpha}(\dd{r}))_{\gamma}=r_{\gamma},
\end{align}
for $\gamma=*$ or $\gamma\in\Omega\setminus\left(\{\alpha,\alpha0,\alpha1\}\cup\{\alpha01\beta,\alpha10\beta:\beta\in\{0,1\}^*\}\right)$.
\end{lemma}
{\bf Proof:} Assume $x\in \mathbb{F}_q$. Let $\theta_{x,\alpha}$ be the map
over $L(\Omega)\cup R(\Omega)$ satisfying (\ref{400}) and, if
$\alpha0\in \Omega$,
\begin{align*}
    (\theta_{x,\alpha}([l]))_{\alpha0}=l_{\alpha0},\ (\theta_{x,\alpha}(\dd{r}))_{\alpha0}=r_{\alpha0}+xr_*,
\end{align*}
if $\alpha1\in \Omega$,
\begin{align*}
    (\theta_{x,\alpha}([l]))_{\alpha1}=l_{\alpha1}-xl_*,\ (\theta_{x,\alpha}(\dd{r}))_{\alpha1}=r_{\alpha1},
\end{align*}
if $\alpha01\beta\in \Omega$,
\begin{align*}
    (\theta_{x,\alpha}([l]))_{\alpha01\beta}=l_{\alpha01\beta}+xl_{\beta},\  (\theta_{x,\alpha}(\dd{r}))_{\alpha01\beta}=r_{\alpha01\beta}+xr_{\beta},
\end{align*}
if $\alpha10\beta\in \Omega$,
\begin{align*}
    (\theta_{x,\alpha}([l]))_{\alpha10\beta}=l_{\alpha10\beta}-xl_{\beta},\ (\theta_{x,\alpha}(\dd{r}))_{\alpha10\beta}=r_{\alpha10\beta}-xr_{\beta},
\end{align*}
and (\ref{401}) for $\gamma=*$ or any other sequence $\gamma$ in
$\Omega$, namely,
\begin{align*}
\gamma\not\in\{\alpha,\alpha0,\alpha1\}\cup\{\alpha01\beta,\alpha10\beta:\beta\in\{0,1\}^*\}.
\end{align*}
Clearly, $\theta_{x,\alpha}$ is well-defined over $L(\Omega)\cup
R(\Omega)$ and a bijective map over $L(\Omega)$ as well as over
$R(\Omega)$. One can also check easily that
$(\theta_{x,\alpha}([l]),\theta_{x,\alpha}(\dd{r}))\in
E(\Gamma(\Omega,q))$ if and only if $([l],\dd{r})\in
E(\Gamma(\Omega,q))$. Hence, the map $\theta_{x,\alpha}$ is indeed
the desired automorphism of $\Gamma(\Omega,q)$.

The proof is completed.
\eproof
\medskip

Let $\mathbb{S}(\Omega)=\cup_{\alpha\in\Omega}S_{\Omega}(\alpha)$. For
any sequence $\alpha\in\{0,1\}^*$ and $i\geq 0$, let $\alpha^i$ denote the
binary sequence defined by
\begin{align*}
    \alpha^i=\left\{
    \begin{array}{ll}
    \underbrace{\alpha\alpha\cdots \alpha}_i\ ,&\text{ if }i>0,\\
    \eta,&\text{ if }i=0.
    \end{array}\right.
\end{align*}

\begin{theorem}
\label{thm:-002}
Assume that $\mathbb{S}(\Omega)\subset\Omega$.
Then, for any two adjacent vertices $[l]\in L(\Omega)$ and $\dd{r}\in
R(\Omega)$,
\begin{enumerate}
    \item there is an automorphism $\theta_0$ of $\Gamma(\Omega,q)$ such that $\theta_0([l])=[0]_{l_*}$ and
\begin{align}
    &(\theta_0(\dd{r}))_*=r_*, \label{703}\\
    &(\theta_0(\dd{r}))_{1^i}=l_*^{i+1}r_*, \text{ if }1^i\in\Omega, i\geq 0,\label{704}\\
    &(\theta_0(\dd{r}))_{\alpha}=0,\text{ for all }\alpha\in\Omega\setminus\{1^i:i\geq 0\},\label{705}
\end{align}
    \item there is an automorphism $\theta_1$ of $\Gamma(\Omega,q)$ such that $\theta_1(\dd{r})=\dd{0}_{r_*}$ and
\begin{align}
    &(\theta_1([l]))_*=l_*, \\
    &(\theta_1([l]))_{0^i}=r_*^{i+1}l_*, \text{ if }0^i\in\Omega, i\geq 0,\\
    &(\theta_1([l]))_{\alpha}=0,\text{ for all }\alpha\in\Omega\setminus\{0^i:i\geq 0\}.
\end{align}
\end{enumerate}
\end{theorem}
{\bf Proof:} We only prove the first result. The other result is
true by symmetry.

Since, for any $\alpha\in\Omega$, the
automorphism $\theta_{x,\alpha}$ given in Lemma~\ref{lem:0002}
fixes those entries of $[l]\in L(\Omega)$ and $\dd{r}\in
R(\Omega)$ whose indices are in $\{*\}\cup\{\beta\in \Omega:
|\beta|\leq |\alpha|,\beta\neq \alpha\}$, we see that there is an
automorphism $\theta_0$ of $\Gamma(\Omega,q)$ satisfying
$\theta_0([l])=[0]_{l_*}$ and (\ref{703}). Furthermore, from
$(\theta_0([l]),\theta_0(\dd{r}))\in E(\Gamma(\Omega,q))$ and
(\ref{000}) to (\ref{002}), we see that (\ref{704}) and
(\ref{705}) are also true. \eproof
\medskip

For any sequence $\alpha\in\{0,1\}^*$, let $H_0(\alpha)$ denote the set
of its subsequences obtained by deleting a symbol 0 either from
the first position or from any two consecutive 0's in a few
iterative stages. For example,
\begin{align*}
    H_0(0^i)&=W_{i+1}=\{0^i,0^{i-1},\ldots,\eta\},\ i\geq 0,\\
    H_0(0^21^20^3)&=\{0^21^20^3, 01^20^3, 1^20^3, 0^21^20^2,01^20^2, 1^20^2, 0^21^20, 01^20, 1^20\}.
\end{align*}
Furthermore, we define a set $T_0(\alpha)$ as the following
\begin{align*}
    &T_0(\alpha)=\left\{
    \begin{array}{ll}
    \{*\},&\text{if }\alpha=0^i,\ i\geq 0,\\
    \emptyset,&\text{if }\alpha=\beta 1,\\
H_0(\beta1),&\text{if }\alpha=\beta10^i,\ i>0.
    \end{array}
    \right.
\end{align*}
Clearly,
\begin{align}
\label{+000}
    T_0(\alpha)\cap H_0(\alpha)=\emptyset.
\end{align}

\begin{lemma}
\label{lem:0003} For any binary sequence $\alpha$,
\begin{align}
\label{-203}
    H_0(\alpha 1)&=\{\beta 1: \beta\in H_0(\alpha)\},\\
\label{-202}
    H_0(\alpha 0)&=\{\beta 0: \beta\in H_0(\alpha)\cup T_0(\alpha)\},\\
\label{-301}
    T_0(\alpha0)&=\left(H_0(\alpha)\cup T_0(\alpha)\right)\setminus H_0(\alpha 0).
\end{align}
\end{lemma}
{\bf Proof:} From the definitions, one can get (\ref{-203}),
(\ref{-202}) and
\begin{align}
\label{1-205}
    T_0(\alpha 0)&=\left\{
    \begin{array}{ll}
    H_0(\alpha),&\text{if }\alpha=\beta1,\\
    T_0(\alpha),&\text{if }\alpha=\eta\text{ or }\beta0,
    \end{array}\right.
\end{align}
immediately. If $\alpha=\beta1$, we have $T_0(\alpha)=\emptyset$. If
$\alpha=\eta$ or $\beta0$, we have $H_0(\alpha)\subset H_0(\alpha 0)$.
Hence, from (\ref{1-205}) we have
\begin{align*}
    &H_0(\alpha 0)\cup T_0(\alpha 0)=H_0(\alpha 0)\cup H_0(\alpha)\cup T_0(\alpha),
\end{align*}
and thus (\ref{-301}) follows from $H_0(\alpha 0)\cap
T_0(\alpha0)=\emptyset$. \eproof
\medskip

Let $\mathbb{H}_0(\Omega)=\bigcup_{\alpha\in\Omega}H_0(\alpha)$.
Similarly, we define $\mathbb{H}_1(\Omega)$ as the set of binary sequences obtained from those in $\Omega$
by deleting a symbol 1 from the first position or two consecutive 1's in a few iterative stages.

\begin{lemma}
\label{lem:-001}
\begin{enumerate}
  \item If $\Omega$ contains all of the sequences in $\mathbb{H}_0(\Omega)$, then, for any $x\in \mathbb{F}_q$, there is an automorphism $\phi$ of
$\Gamma(\Omega,q)$ such that
\begin{align*}
(\phi([l]))_*&=l_*,\text{ for }[l]\in L(\Omega),\\
(\phi(\dd{r}))_*&=r_*+x,\text{ for }\dd{r}\in R(\Omega).
\end{align*}
  \item If $\Omega$ contains all of the sequences in $\mathbb{H}_1(\Omega)$, then, for any $x\in \mathbb{F}_q$, there is an automorphism $\psi$ of
$\Gamma(\Omega,q)$ such that
\begin{align*}
(\psi([l]))_*&=l_*+x,\text{ for }[l]\in L(\Omega),\\
(\psi(\dd{r}))_*&=r_*,\text{ for }\dd{r}\in R(\Omega).
\end{align*}
\end{enumerate}

\end{lemma}
{\bf Proof:} We only prove the first result. The other result is
true by symmetry.

Suppose $\mathbb{H}_0(\Omega)\subset \Omega$. Since $\Omega$ satisfies the condition C1, we see easily $H_0(\alpha)\cup T_0(\alpha)\subset \Omega\cup\{*\}$
for any $\alpha\in\Omega$.
Let $[l]\in L(\Omega)$, $\dd{r}\in
R(\Omega)$ be two adjacent vertices of $\Gamma(\Omega,q)$ and $x$ an
element in $\mathbb{F}_q$.

At first, we define $f_{x,\eta}(\eta)=1$,
$f_{x,\eta}(*)=x$ and $g_{x,\eta}(\eta)=1$. Then,
$f_{x,\eta}(\cdot)$ and $g_{x,\eta}(\cdot)$ have been well defined
over $H_0(\eta)\cup T_0(\eta)$ and $H_0(\eta)$, respectively, and
\begin{align}
    &\sum_{\gamma\in H_0(\eta)\cup T_0(\eta)}f_{x,\eta}(\gamma)l_{\gamma}+\sum_{\gamma\in
    H_0(\eta)}g_{x,\eta}(\gamma)r_{\gamma}=l_{\eta}+xl_*+r_{\eta}
    =l_*(r_*+x).
    \label{600}
\end{align}

Now we assume that $f_{x,\alpha}(\cdot)$ and $g_{x,\alpha}(\cdot)$ have
been well defined for some $\alpha\in\{0,1\}^*$ over $H_0(\alpha)\cup
T_0(\alpha)$ and over $H_0(\alpha)$, respectively.

For $\beta\in T_0(\alpha0)$, let
$f_{x,\alpha0}(\beta)=xf_{x,\alpha}(\beta)$. For $\beta\in
H_0(\alpha)\cup T_0(\alpha)$, let
$g_{x,\alpha0}(\beta0)=f_{x,\alpha}(\beta)$ and
\begin{align*}
    f_{x,\alpha0}(\beta0)&=\left\{
    \begin{array}{ll}
    f_{x,\alpha}(\beta)+xf_{x,\alpha}(\beta0),&\text{if }\beta0\in H_0(\alpha)\cup T_0(\alpha),\\
    f_{x,\alpha}(\beta),&\text{otherwise}.
    \end{array}\right.
\end{align*}
Then, from (\ref{-202}) we see that $f_{x,\alpha0}(\cdot)$ and
$g_{x,\alpha0}(\cdot)$ have been well defined over $H_0(\alpha0)\cup
T_0(\alpha0)$ and $H_0(\alpha0)$, respectively. Furthermore, if
$\alpha0\in\Omega$, from (\ref{-301}) we have
\begin{align}
    &\sum_{\gamma\in H_0(\alpha0)\cup T_0(\alpha0)}f_{x,\alpha 0}(\gamma)l_{\gamma}+\sum_{\gamma\in H_0(\alpha0)}g_{x,\alpha 0}(\gamma)r_{\gamma}\nonumber\\
    =&\sum_{\beta\in H_0(\alpha)\cup T_0(\alpha)}f_{x,\alpha0}(\beta0)l_{\beta0}
    +\sum_{\beta\in T_0(\alpha0)}f_{x,\alpha0}(\beta)l_{\beta}+\sum_{\beta\in H_0(\alpha)\cup T_0(\alpha)}g_{x,\alpha0}(\beta0)r_{\beta0}\nonumber\\
    =&\sum_{\beta\in H_0(\alpha)\cup T_0(\alpha)}f_{x,\alpha}(\beta)l_{\beta0}
    +\sum_{\beta0\in (H_0(\alpha)\cup T_0(\alpha))\cap H_0(\alpha0)}xf_{x,\alpha}(\beta0)l_{\beta0}\nonumber\\
    &\hspace{2.5cm}+\sum_{\beta\in (H_0(\alpha)\cup T_0(\alpha))\setminus H_0(\alpha0)}xf_{x,\alpha}(\beta)l_{\beta}+
    \sum_{\beta\in H_0(\alpha)\cup T_0(\alpha)}f_{x,\alpha}(\beta)r_{\beta0}\nonumber\\
    =&\sum_{\beta\in H_0(\alpha)\cup T_0(\alpha)}f_{x,\alpha}(\beta)(l_{\beta0}+r_{\beta0}+xl_{\beta})\nonumber\\
    =&(r_*+x)\sum_{\beta\in H_0(\alpha)\cup T_0(\alpha)}f_{x,\alpha}(\beta)l_{\beta}.
    \label{601}
\end{align}

For $\beta\in H_0(\alpha)$, let
$f_{x,\alpha1}(\beta1)=g_{x,\alpha1}(\beta1)=g_{x,\alpha}(\beta)$.
Then, from (\ref{-203}) and $T_0(\alpha1)=\emptyset$, we see that
$f_{x,\alpha1}(\cdot)$ and $g_{x,\alpha1}(\cdot)$ have been well
defined over $H_0(\alpha1)\cup T_0(\alpha1)$ and over $H_0(\alpha1)$,
respectively, and, if $\alpha1\in\Omega$,
\begin{align}
    &\sum_{\gamma\in H_0(\alpha1)\cup T_0(\alpha1)}f_{x,\alpha 1}(\gamma)l_{\gamma}+\sum_{\gamma\in H_0(\alpha1)}g_{x,\alpha
    1}(\gamma)r_{\gamma}\nonumber\\
    =&\sum_{\beta\in H_0(\alpha)}g_{x,\alpha}(\beta)(l_{\beta1}+r_{\beta1})
    =l_*\sum_{\beta\in H_0(\alpha)}g_{x,\alpha}(\beta)r_{\beta}.
    \label{602}
\end{align}

Then, $f_{x,\alpha}(\cdot)$ and $g_{x,\alpha}(\cdot)$ have been
well defined over $H_0(\alpha)\cup T_0(\alpha)$ and over $H_0(\alpha)$,
respectively, for all $\alpha\in\{0,1\}^*$.

Let $\phi([l])$ denote
the vector in $L(\Omega)$ defined by $(\phi([l]))_*=l_*$ and
\begin{align}
    (\phi([l]))_{\alpha}=\sum_{\beta\in H_0(\alpha)\cup T_0(\alpha)}f_{x,\alpha}(\beta)l_{\beta}, \text{ for }\alpha\in\Omega.
\end{align}
Let $\phi(\dd{r})$ denote the vector in $R(\Omega)$ defined by
$(\phi(\dd{r}))_*=r_*+x$ and
\begin{align}
    (\phi(\dd{r}))_{\alpha}=\sum_{\beta\in H_0(\alpha)}g_{x,\alpha}(\beta)r_{\beta}, \text{ for }\alpha\in\Omega.
\end{align}
Therefore, from (\ref{600}) to (\ref{602}), we see that
$\phi([l])\in L(\Omega)$ and $\phi(\dd{r})\in R(\Omega)$ are
adjacent in $\Gamma(\Omega,q)$. Clearly, for any $\alpha\in
\Omega$, we have $f_{x,\alpha}(\alpha)=g_{x,\alpha}(\alpha)=1$ and
that any sequence in $(H_0(\alpha)\cup
T_0(\alpha))\setminus\{\alpha\}$ is shorter than $\alpha$. Hence,
$\phi([l])$ and $\phi(\dd{r})$ are bijective maps in $L(\Omega)$
and $R(\Omega)$, respectively. Thus, $\phi$ is the desired
automorphism of $\Gamma(\Omega,q)$. \eproof
\medskip

We note that the automorphism $\phi$ in Lemma~\ref{lem:-001} is partially given by, for $k\geq 0$,
\begin{align*}
    &\left\{\begin{array}{l}
    (\phi([l]))_{(10)^k}=l_{(10)^k}+xl_{*(01)^k},\\
    (\phi(\dd{r})_{(10)^k}=r_{(10)^k},
    \end{array}
    \right. \text{ if }(10)^k\in\Omega,\\
    &\left\{\begin{array}{l}
    (\phi([l]))_{0(10)^k}=l_{0(10)^k}+xl_{(10)^k}+xl_{(01)^k}+x^2l_{*(01)^k},\\
    (\phi(\dd{r})_{0(10)^k}=r_{0(10)^k}+xr_{(10)^k},
    \end{array}
    \right. \text{ if }0(10)^k\in\Omega,\\
    &\left\{\begin{array}{l}
    (\phi([l]))_{1(01)^k}=l_{1(01)^k},\\
    (\phi(\dd{r})_{1(01)^k}=r_{1(01)^k},
    \end{array}
    \right. \text{ if }1(01)^k\in\Omega,\\
        &\left\{\begin{array}{l}
    (\phi([l]))_{(01)^{k+1}}=l_{(01)^{k+1}}+xl_{1(01)^k},\\
    (\phi(\dd{r})_{(01)^{k+1}}=r_{(01)^{k+1}}+xr_{1(01)^k},
    \end{array}
    \right. \text{ if }(01)^{k+1}\in\Omega,
\end{align*}
where $*(01)^k$ denotes the symbol $*$ if $k=0$, and the sequence $1(01)^{k-1}$ otherwise.

\begin{theorem}
\label{thm:401}
Assume that $\mathbb{S}(\Omega)\subset\Omega$.
\begin{enumerate}
  \item If $\Omega$ contains all of the sequences in $\mathbb{H}_0(\Omega)$, then, for any $\dd{r}\in R(\Omega)$, there is an automorphism $\pi$ of
$\Gamma(\Omega,q)$ such that $\pi(\dd{r})=\dd{0}_0$.
 \item If $\Omega$ contains all of the sequences in $\mathbb{H}_1(\Omega)$, then, for any $[l]\in L(\Omega)$, there is an automorphism $\pi$ of
$\Gamma(\Omega,q)$ such that $\pi([l])=[0]_0$.
\end{enumerate}
\end{theorem}
{\bf Proof:} We only prove the first result. The other result is
true by symmetry.

Suppose $\mathbb{H}_0(\Omega)\subset \Omega$. Let $\dd{r}\in$ be an arbitrary vertex in $R(\Omega)$.
According to Lemma~\ref{lem:-001}, there is an automorphism $\phi$ of $\Gamma(\Omega,q)$
such that $(\phi(\dd{r}))_*=0$.
Hence, $\pi=\theta_1\phi$ is the
desired automorphism of $\Gamma(\Omega,q)$, where $\theta_1$ is
the automorphism given in Theorem~\ref{thm:-002} for the vertex $\phi(\dd{r})$. \eproof
\medskip

\begin{theorem}
\label{thm:-001}If $\Omega$ contains all of the sequences in $\mathbb{H}_0(\Omega)\cup\mathbb{H}_1(\Omega)$,
then, the bipartite graph $\Gamma(\Omega,q)$ is edge-transitive,
or equivalently, for any pair of adjacent vertices $[l]\in
L(\Omega)$ and $\dd{r}\in R(\Omega)$ there is an automorphism
$\pi$ of $\Gamma(\Omega,q)$ such that $\pi([l])=[0]_0$ and
$\pi(\dd{r})=\dd{0}_0$.
\end{theorem}
{\bf Proof:}
Suppose that $\mathbb{H}_0(\Omega)\cup\mathbb{H}_1(\Omega)\subset\Omega$.
Since $\mathbb{H}_0(\Omega)\cup\mathbb{H}_1(\Omega)$ contains
all of the sequences obtained from those in $\Omega$ by deleting the first bit, we see that $S_{\Omega}(\alpha)\subset\Omega$ is valid for all $\alpha\in\Omega$.

Assume that $[l]\in L(\Omega)$ and $\dd{r}\in
R(\Omega)$ are adjacent in $\Gamma(\Omega,q)$.
According to Lemma~\ref{lem:-001}, there are automorphisms $\phi$, $\psi$ of $\Gamma(\Omega,q)$
such that
\begin{align*}
    (\phi(\dd{r}))_*&=0,\\
    (\phi([l]))_*&=l_*,\\
    (\psi\phi([l]))_*&=0,\\
    (\psi\phi(\dd{r}))_*&=(\phi(\dd{r}))_*=0.
\end{align*}
Hence, $\pi=\theta_0\psi\phi$ is the
desired automorphism of $\Gamma(\Omega,q)$, where $\theta_0$ is
the automorphism given in Theorem~\ref{thm:-002}
for the adjacent vertices $\psi\phi([l])$ and $\psi\phi(\dd{r})$. \eproof
\medskip

\begin{corollary}
For $k\geq 2$, the bipartite graph $D(k,q)$ is edge-transitive.
\end{corollary}
{\bf Proof:} For $k\geq 2$, it is very easy to check that $\mathbb{H}_0(\Omega)\cup\mathbb{H}_1(\Omega)\subset\Omega$ is valid for $\Omega=U_k$, where $U_k$ is the set consisting of the first $k-1$ sequences in the set $U$ defined by (\ref{0009}).
Then, according to Theorem~\ref{thm:-001}
the bipartite graph $\Gamma(U_k,q)$, which is equivalent to $D(k,q)$, is edge-transitive. \eproof
\medskip

\noindent {\bf Example 1}: For any finite set $\Lambda$ of binary
sequences, let $\Phi(\Lambda)$ denote the smallest set such that
$\Lambda\subset\Phi(\Lambda)$ and $\mathbb{H}_0(\Phi(\Lambda))\cup\mathbb{H}_1(\Phi(\Lambda))\subset\Phi(\Lambda)$. Let
$\Lambda_0$ be a finite set of binary sequence such that, for any
$\alpha\in \Lambda_0$, any sequence obtained from $\alpha$ by
deleting  a bit either from the first position, or from two
consecutive 0's, or from two consecutive 1's, is still in
$\Lambda_0$. Let $\Lambda_1$ be a finite set of binary sequences
such that, for any $\alpha\in \Lambda_1$, any sequence obtained
from $\alpha$ by deleting a bit either from the last position, or
from two consecutive 0's, or from two consecutive 1's, is still
in $\Lambda_1$. If there is a symbol $a\in \{0,1\}$ such that $a\not\in \Lambda_0$ and $\bar{a}\not\in\Lambda_1$,
where $\bar{a}$ is the symbol in $\{0,1\}$ other than $a$, then the set
\begin{align}\label{701}
    \Omega_1=\Phi(\Lambda_0)\cup \Phi(\Lambda_1)\cup\{\alpha\beta:\alpha\in \Lambda_0,\beta\in\Lambda_1\}
\end{align}
satisfies $\mathbb{H}_0(\Omega_1)\cup\mathbb{H}_1(\Omega_1)\subset\Omega_1$, and thus from Theorem~\ref{thm:-001}
the bipartite graph $\Gamma(\Omega_1,q)$ is edge-transitive.

For example, if we take $\Lambda_0=\{\eta,0,10,0^2\}$ and
$\Lambda_1=\{\eta,1,10,101,(10)^2,(10)^21\}$, then $\Phi(\Lambda_0)=\{\eta,0,1,0^2,10\}$, $\Phi(\Lambda_1)=\{\eta,0,1,01,10,010,101,(01)^2,(10)^2,(10)^21\}$ and the set
given by (\ref{701}), denoted $\Omega_2$, consists of the following sequences
\begin{align*}
    \eta,&1,0,10,0^2,01,101,0^21,010,(10)^2,0^210,(01)^2,(10)^21,\\
    &0(01)^2,(01)^20,(10)^3,0^2(10)^2,(01)^3,1(01)^3,0(01)^3.
\end{align*}
If the above sequences are mapped into integers $1,2,\ldots,20$ in
order and the symbol $*$ is mapped to 0, then $[l]$
and $\dd{r}$ are adjacent in the edge-transitive bipartite graph
$\Gamma(\Omega_2,q)$ if and only if
\begin{align*}
    l_1+r_1=l_0r_0,\\
    l_2+r_2=l_0r_1,\\
    l_3+r_3=r_0l_1,\\
    l_4+r_4=r_0l_2,\\
    l_5+r_5=r_0l_3,
\end{align*}
and
\begin{align*}
    l_i+r_i=\left\{
    \begin{array}{ll}
    l_0r_{i-3},&\text{if }i\equiv 0\text{ or }1\text{ or }2(\bmod 6),\\
    r_0l_{i-3},&\text{if }i\equiv 3\text{ or }4\text{ or }5(\bmod 6),
    \end{array}\right.
\end{align*}
for $6\leq i\leq 20$.\eproof

\section{Connectivity of $\Gamma(\Omega,q)$}\label{sec000}
In this section, we mainly discuss the connectivity of the graph $\Gamma(\Omega,q)$.
For the components of $\Gamma(\Omega,q)$, a lower bound for their amount is given by showing some of their invariants at first, and their behaviors
under some projections naturally defined are investigated then.

We will assume in general that $[l]\in L(\Omega)$ and $\dd{r}\in
R(\Omega)$ are two adjacent vertices of $\Gamma(\Omega,q)$ without specification in this section.

For any $x\in \mathbb{F}_q$, let $x^0$ be the multiplicative unit of $\mathbb{F}_q$.

\begin{lemma}
\label{lem:d000}
Suppose $\alpha,\beta\in\{0,1\}^*\cup\{*\}$. For any nonnegative integer $s$, we have
\begin{enumerate}
  \item If $\{\alpha 10^s, \beta 10^s\}\subset\Omega$, then
\begin{align}
\label{y001} r_{\beta}l_{\alpha10^s}-r_{\alpha}l_{\beta10^s}
=r_{\alpha}\sum_{t=0}^sr_*^{s-t}r_{\beta10^t}
-r_{\beta}\sum_{t=0}^sr_*^{s-t}r_{\alpha10^t}.
\end{align}
  \item If $\{\alpha 01^s, \beta 01^s\}\subset\Omega$, then
\begin{align}
\label{y002} l_{\beta}r_{\alpha01^s}-l_{\alpha}r_{\beta01^s}
=l_{\alpha}\sum_{t=0}^sl_*^{s-t}l_{\beta01^t}
-l_{\beta}\sum_{t=0}^sl_*^{s-t}l_{\alpha01^t}.
\end{align}
\end{enumerate}

\end{lemma}
{\bf Proof:} We only prove (\ref{y001}). A proof for (\ref{y002})
can be given similarly.

If $\{\alpha 1, \beta 1\}\subset\Omega$, from (\ref{000}) to
(\ref{002}) we see
\begin{align*}
r_{\beta}l_{\alpha1}-r_{\alpha}l_{\beta1}
=r_{\beta}(l_*r_{\alpha}-r_{\alpha1})-r_{\alpha}(l_*r_{\beta}-r_{\beta1})
=r_{\alpha}r_{\beta1}-r_{\beta}r_{\alpha1}.
\end{align*}
Hence, (\ref{y001}) is valid for $s=0$.

Now we assume that (\ref{y001}) is valid for some $s$ with $s\geq
0$. If $\{\alpha 10^{s+1}, \beta 10^{s+1}\}\subset\Omega$,
\begin{align*}
&r_{\beta}l_{\alpha10^{s+1}}-r_{\alpha}l_{\beta10^{s+1}}\\
=&r_{\beta}(r_*l_{\alpha10^s}-r_{\alpha10^{s+1}})-r_{\alpha}(r_*l_{\beta10^s}-r_{\beta10^{s+1}})\\
=&(r_{\alpha}r_{\beta10^{s+1}}-r_{\beta}r_{\alpha10^{s+1}})+r_*(r_{\beta}l_{\alpha10^s}-r_{\alpha}l_{\beta10^s})\\
=&(r_{\alpha}r_{\beta10^{s+1}}-r_{\beta}r_{\alpha10^{s+1}})+
r_*\left(r_{\alpha}\sum_{t=0}^sr_*^{s-t}r_{\beta10^t}
-r_{\beta}\sum_{t=0}^sr_*^{s-t}r_{\alpha10^t}\right)\\
=&r_{\alpha}\sum_{t=0}^{s+1}r_*^{s+1-t}r_{\beta10^t}
-r_{\beta}\sum_{t=0}^{s+1}r_*^{s+1-t}r_{\alpha10^t}.
\end{align*}
Hence, (\ref{y001}) is valid for $s+1$.

The proof is completed.\eproof
\medskip

\begin{lemma}
\label{lem:d001}
Suppose $\alpha\in\{0,1\}^*\cup\{*\}$.
\begin{enumerate}
  \item If $\alpha0^q$ is a sequence in $\Omega$, then
\begin{align}
l_{\alpha0^{q}}-l_{\alpha0}=r_{\alpha0}-(r_{\alpha0^q}+r_*r_{\alpha0^{q-1}}+\cdots+r_*^{q-1}r_{\alpha0})\label{y10}
\end{align}
is an invariant for each component of the graph $\Gamma(\Omega,q)$.
  \item If $\alpha1^q$ is a sequence in $\Omega$, then
\begin{align}
r_{\alpha1^{q}}-r_{\alpha1}=l_{\alpha1}-(l_{\alpha1^q}+l_*l_{\alpha1^{q-1}}+\cdots+l_*^{q-1}l_{\alpha1})\label{y11}
\end{align}
is an invariant for each component of the graph $\Gamma(\Omega,q)$.
\end{enumerate}
\end{lemma}
{\bf Proof}:
We only prove (\ref{y10}). A proof for (\ref{y11})
can be given similarly. Clearly, we have
\begin{align*}
l_{\alpha0^{q}}
&=r_*l_{\alpha 0^{q-1}}-r_{\alpha0^q}\\
&=r_*^2l_{\alpha 0^{q-2}}-(r_{\alpha0^q}+r_*r_{\alpha0^{q-1}})\\
&=r_*^ql_{\alpha}-(r_{\alpha0^q}+r_*r_{\alpha0^{q-1}}+\cdots+r_*^{q-1}r_{\alpha0})\\
&=r_*l_{\alpha}-(r_{\alpha0^q}+r_*r_{\alpha0^{q-1}}+\cdots+r_*^{q-1}r_{\alpha0})\\
&=l_{\alpha0}+r_{\alpha0}-(r_{\alpha0^q}+r_*r_{\alpha0^{q-1}}+\cdots+r_*^{q-1}r_{\alpha0}),
\end{align*}
thus (\ref{y10}) follows.\eproof
\medskip

From Lemmas~\ref{lem:d000} and~\ref{lem:d001}, one can deduce the following corollary easily.
\begin{corollary}
\label{cor:y300}
Suppose $\alpha,\beta\in\{0,1\}^*$. Let $s$ and $s'$ be positive integers with $(q-1)|s-s'$.
\begin{enumerate}
  \item If $\{\alpha 10^s, \beta 10^{s'}\}\subset\Omega$, then
\begin{align}\label{cor001}
    &r_{\beta}l_{\alpha 10^{s}}-r_{\alpha}l_{\beta 10^{s'}}
    =r_{\alpha}\sum_{t=0}^{s'}r_*^{s'-t}r_{\beta 10^{t}}-r_{\beta}\sum_{t=0}^{s}r_*^{s-t}r_{\alpha 10^{t}}.
\end{align}
  \item If $\{\alpha 01^s, \beta 01^{s'} \}\subset\Omega$, then
\begin{align}\label{cor002}
    &l_{\beta}r_{\alpha 01^{s}}-l_{\alpha}r_{\beta 01^{s'}}
    =l_{\alpha}\sum_{t=0}^{s'}l_*^{s'-t}l_{\beta 01^{t}}-l_{\beta}\sum_{t=0}^{s}l_*^{s-t}l_{\alpha 01^{t}}.
\end{align}
\end{enumerate}
Furthermore, if either $\alpha$ or $\beta$ is $*$, then (\ref{cor001}) and (\ref{cor002}) are valid for
nonnegative integers $s$, $s'$ with $(q-1)|s-s'$.
\end{corollary}
\medskip

For symbol $a\in\{0,1\}$ and sequence $\boldsymbol{s}=(s_1,s_2,\ldots)$ of
nonnegative integers, let $\mu_{a,0}(\boldsymbol{s})=*$ and, for
$i\geq 1$,
\begin{align}
\mu_{a,i}(\boldsymbol{s})=\left\{
\begin{array}{ll}
\mu_{a,i-1}(\boldsymbol{s})10^{s_{i}},&\text{if }i+a\text{ is odd},\\
\mu_{a,i-1}(\boldsymbol{s})01^{s_{i}},&\text{if }i+a\text{ is even},
\end{array}\right.
\end{align}
where $a$ is also treated as an integer.
Clearly, for any nonempty sequence $\alpha\in \{0,1\}^*$,
there uniquely exist a symbol $a\in\{0,1\}$, a positive integer $n$ and nonnegative integers $s_1,s_2,\ldots,s_n$
such that
\begin{align}\label{krs01}
\alpha=\mu_{a,n+1}(s_1,s_2,\ldots,s_n,0)=\mu_{a,n}(s_1,\ldots,s_{n-1},s_n+1).
\end{align}
For such $\alpha$, we define $\kappa(\alpha)=n$ and
\begin{align}\label{krs02}
    \zeta_{\text{L}}(\alpha,[l])&=\left\{
    \begin{array}{ll}
    \sum_{t=0}^{s_n+1}l_*^{s_n+1-t}l_{\mu_{a,n-1}(s_1,\ldots,s_{n-1})01^t},&\text{ if } n+a \text{ is even},\\
    l_{\alpha},&\text{ if } n+a \text{ is odd},
    \end{array}
    \right.\\
    \zeta_{\text{R}}(\alpha,\dd{r})&=\left\{
    \begin{array}{ll}
    \sum_{t=0}^{s_n+1}r_*^{s_n+1-t}r_{\mu_{a,n-1}(s_1,\ldots,s_{n-1})10^t},&\text{ if } n+a \text{ is odd},\\
    r_{\alpha},&\text{ if } n+a \text{ is even},
    \end{array}
    \right.
\end{align}
where $a$ is still treated as an integer.

For positive integer $n$, let $\Delta_n(\Omega)=\{\beta\in\Omega^*:\kappa(\beta)<n\}\cup\{*,\eta\}$.
For polynomials $P, Q$ of the entries of the vertices $[l], \dd{r}$,
we write
\begin{align}
    P\stackrel{n}{=}Q
\end{align}
if there are polynomials $f$
and $g$ of the variables in $\{l_{\gamma}:\gamma\in\Delta_n(\Omega)\}$ and in $\{r_{\gamma}:\gamma\in\Delta_n(\Omega)\}$, respectively, such that
\begin{align}
    P-Q=f+g.
\end{align}

\begin{theorem}
\label{thm:m01}
Let $n$ be a positive integer.
Assume that $\boldsymbol{s}=(s_1,\ldots,s_n,0)$, $\boldsymbol{s}'=(s'_1,\ldots,s'_n,0)$ and $\boldsymbol{s}''=(s''_1,\ldots,s''_n,0)$ are sequences of nonnegative integers
such that
\begin{align}
    &s_i\equiv s_i'\equiv s_{n+1-i}''\bmod (q-1), \text{ for }i=1,2,\ldots,n,\\
    &s_i=s_i'=s_{n+1-i}'',\text{ if }\min\{s_i, s_i', s_{n+1-i}''\}=0\text{ and } 1<i<n.
\end{align}
Suppose $a, b$ are symbols in $\{0,1\}$ such that $a=b$ if and only if $n$ is odd. Let $\alpha=\mu_{a,n+1}(\boldsymbol{s})$, $\beta=\mu_{a,n+1}(\boldsymbol{s}')$ or $\mu_{b,n+1}(\boldsymbol{s}'')$.
If $\{\alpha,\beta,\mu_{b,n-m}(\boldsymbol{s}'')\}\subset\Omega$,
where $m$ is the smallest integer such that $s_m\neq s_m'$, then for each edge
$([l],\dd{r})\in E(\Gamma(\Omega,q))$ we have
\begin{align}
\label{gg00}
    \zeta_{\text{L}}(\alpha,[l])+\zeta_{\text{R}}(\alpha,\dd{r})
    \stackrel{n}{=}
    \zeta_{\text{L}}(\beta,[l])+\zeta_{\text{R}}(\beta,\dd{r}).
\end{align}
\end{theorem}
{\bf Proof}:
Without loss of generality we assume that $a=1$ and $2|n$.
Then, $b=0$ and
\begin{align*}
\alpha=\mu_{1,n+1}(\boldsymbol{s})=1^{s_1+1}0^{s_2+1}\ldots1^{s_{n-1}+1}0^{s_{n}+1}.
\end{align*}
If $\beta=\mu_{0,n+1}(\boldsymbol{s}'')=0^{s_1''+1}1^{s_2''+1}\cdots0^{s_{n-1}''+1}1^{s_{n}''+1}$,
from Corollary~\ref{cor:y300} we see
\begin{align*}
&\zeta_{\text{L}}(\alpha,[l])+\zeta_{\text{R}}(\alpha,\dd{r})\\
=&l_{\alpha}+r_{\alpha}+
\sum_{t=0}^{s_n}r_*^{s_n+1-t}r_{\mu_{1,n-1}(\boldsymbol{s})10^{t}}\\
    =&r_{*}l_{\mu_{1,n-1}(\boldsymbol{s})10^{s_{n}}}+r_*\sum_{t=0}^{s_n}r_*^{s_n-t}r_{\mu_{1,n-1}(\boldsymbol{s})10^{t}}\\
    \stackrel{n}{=}&l_{0^{s''_1}}r_{\mu_{1,n-1}(\boldsymbol{s})}= l_{\mu_{0,1}(\boldsymbol{s}'')}r_{\mu_{1,n-2}(\boldsymbol{s})01^{s_{n-1}}}
    \\
    \stackrel{n}{=} &r_{\mu_{0,1}(\boldsymbol{s}'')01^{s''_2}}l_{\mu_{1,n-2}(\boldsymbol{s})}
    =r_{\mu_{0,2}(\boldsymbol{s}'')}l_{\mu_{1,n-3}(\boldsymbol{s})10^{s_{n-2}}}
    \\
    \vdots\ &\\
    \stackrel{n}{=}
    &l_{\mu_{0,n-2}(\boldsymbol{s}'')10^{s_{n-1}''}}r_{\mu_{1,1}(\boldsymbol{s})}
    =l_{\mu_{0,n-1}(\boldsymbol{s}'')}r_{1^{s_{1}}}\\
    \stackrel{n}{=}&l_{*}r_{\mu_{0,n-1}(\boldsymbol{s''})01^{s_{n}''}}+
    l_*\sum_{t=0}^{s_n''}l_*^{s_n''-t}l_{\mu_{0,n-1}(\boldsymbol{s''})01^{t}}\\
    =&\zeta_{\text{L}}(\beta,[l])+\zeta_{\text{R}}(\beta,\dd{r}).
\end{align*}
If $\beta=\mu_{1,n+1}(\boldsymbol{s}')=1^{s_1'+1}0^{s_2'+1}\cdots1^{s_{n-1}'+1}0^{s_{n}'+1}$,
from Corollary~\ref{cor:y300} we also see
\begin{align*}
&\zeta_{\text{L}}(\alpha,[l])+\zeta_{\text{R}}(\alpha,\dd{r})\\
\stackrel{n}{=} &\left\{
    \begin{array}{ll}
    l_{\mu_{0,n-m}(\boldsymbol{s}'')}r_{\mu_{1,m-1}(\boldsymbol{s})01^{s_m}},&\text{if }m\text{ is odd},\\
    r_{\mu_{0,n-m}(\boldsymbol{s}'')}l_{\mu_{1,m-1}(\boldsymbol{s})10^{s_m}},&\text{if }m\text{ is even},
\end{array}\right.\\
= &\left\{
    \begin{array}{ll}
    l_{\mu_{0,n-m}(\boldsymbol{s}'')}r_{\mu_{1,m-1}(\boldsymbol{s'})01^{s_m}},&\text{if }m\text{ is odd},\\
    r_{\mu_{0,n-m}(\boldsymbol{s}'')}l_{\mu_{1,m-1}(\boldsymbol{s'})10^{s_m}},&\text{if }m\text{ is even},
\end{array}\right.\\
\stackrel{n}{=} &\left\{
    \begin{array}{ll}
    l_{\mu_{0,n-m}(\boldsymbol{s}'')}r_{\mu_{1,m-1}(\boldsymbol{s'})01^{s_m'}},&\text{if }m\text{ is odd},\\
    r_{\mu_{0,n-m}(\boldsymbol{s}'')}l_{\mu_{1,m-1}(\boldsymbol{s'})10^{s_m'}},&\text{if }m\text{ is even},
\end{array}\right.\\
\stackrel{n}{=}&\zeta_{\text{L}}(\beta,[l])+\zeta_{\text{R}}(\beta,\dd{r}),
\end{align*}
where the third equality is deduced according to Lemma~\ref{lem:d001}.
\eproof
\medskip

For the sequences $\alpha$ and $\beta$ satisfying the condition of Theorem~\ref{thm:m01}, we write
$\alpha\bowtie_q\beta$. Clearly, if $\alpha\neq\beta$ and $\alpha\bowtie_q\beta$, then (\ref{gg00}) implies a nontrival invariant
for each component of $\Gamma(\Omega,q)$ of form
\begin{align*}
\zeta_{\text{L}}(\alpha,[l])-\zeta_{\text{L}}(\beta,[l])+f_{\alpha,\beta}
=\zeta_{\text{R}}(\beta,\dd{r})-\zeta_{\text{R}}(\alpha,\dd{r})+g_{\alpha,\beta},
\end{align*}
where $f_{\alpha,\beta}$ and $g_{\alpha,\beta}$ are polynomials of the variables in $\{l_{\gamma}:\gamma\in\Delta_n(\Omega)\}$ and in $\{r_{\gamma}:\gamma\in\Delta_n(\Omega)\}$, respectively, and
$n=\kappa(\alpha)=\kappa(\beta)$.
For example, if $\alpha=\mu_{0,2k+1}(\boldsymbol{0})=(01)^k$ and $\beta=\mu_{1,2k+1}(\boldsymbol{0})=(10)^k$, then the invariant implied in (\ref{gg00}) is
\begin{align*}
    &l_{(01)^k}-l_{(10)^k}+l_*l_{(01)^{k-1}0}+\sum_{i=0}^{k-2}l_{(10)^i1}l_{(01)^{k-2-i}0}-\sum_{i=0}^{k-1}l_{(10)^i}l_{(01)^{k-1-i}}\\
    =&r_{(10)^k}-r_{(01)^k}+r_*r_{(10)^{k-1}1}+\sum_{i=0}^{k-2}r_{(10)^i1}r_{(01)^{k-2-i}0}-\sum_{i=0}^{k-1}r_{(10)^i}r_{(01)^{k-1-i}}.
\end{align*}

Clearly, $\bowtie_q$ defines an equivalent relation in $\Omega^*$
and divides $\Omega^*$ into subsets $\aleph_1,\aleph_2,\ldots,\aleph_{\hbar}$ such that
each of them is an equivalent class defined by $\bowtie_q$, where $\hbar=\hbar(\Omega,q)$ is the number of such equivalent classes. We note that $\aleph_i$ may contain only one sequence.

Now we can give a lower bound for the number of the components of $\Gamma(\Omega,q)$.

\begin{corollary}
The number of components of $\Gamma(\Omega,q)$ is at least $q^{|\Omega^*|-\hbar(\Omega,q)}$.
\end{corollary}
{\bf Proof}:
This corollary follows immediately from (\ref{gg00}) and that, for any subset $Y$ of $\Omega^*$,
\begin{align}\label{r00}
    \{(\zeta_{\text{L}}(\alpha,[l]):\alpha\in Y):[l]\in L(\Omega)\}=\{(\zeta_{\text{L}}(\alpha,\dd{r}):\alpha\in Y):\dd{r}\in R(\Omega)\}
\end{align}
is a vector space over $\mathbb{F}_q$ of dimension $|Y|$.
\eproof\medskip

In the rest of this section, we suppose that $\Omega'$ is a subset of $\Omega$ satisfying the condition C1.
Let $\Pi_{\Omega/\Omega'}$ denote the projection from $\Gamma(\Omega,q)$ to $\Gamma(\Omega',q)$
defined naturally. For a component $C$ of $\Gamma(\Omega,q)$ and vertex $u\in V(\Pi_{\Omega/\Omega'}(C))$, let $\mathcal{V}(u,C)$ denote the set of vertices $v\in V(C)$ with $\Pi_{\Omega/\Omega'}(v)=u$.

\begin{theorem}\label{gg01}
For any component $C$ of $\Gamma(\Omega,q)$, $\Pi_{\Omega/\Omega'}(C)$ is a component of $\Gamma(\Omega',q)$ and
$\Pi_{\Omega/\Omega'}$ is a $t$-to-1 graph homomorphism from $C$ to $\Pi_{\Omega/\Omega'}(C)$ for some $t$ with $1\leq t\leq q^{|\Omega|-|\Omega'|}$.
\end{theorem}
{\bf Proof}:
Let $C'=\Pi_{\Omega/\Omega'}(C)$.
Since for any two adjacent vertices $v, v'$ of $\Gamma(\Omega,q)$, their projections $\Pi_{\Omega/\Omega'}(v), \Pi_{\Omega/\Omega'}(v')$ are adjacent in $\Gamma(\Omega',q)$,
it is clearly that $C'$ is a connected subgraph of $\Gamma(\Omega',q)$.
If $u\in V(C')$ and $u'\in V(\Gamma(\Omega',q))$ are adjacent in $\Gamma(\Omega',q)$, it is clear that, for any $v\in \mathcal{V}(u,C)$,
there is a unique vertex $v'\in V(\Gamma(\Omega,q))$ with $\Pi_{\Omega/\Omega'}(v')=u'$
such that $v'$ and $v$ are adjacent in $\Gamma(\Omega,q)$. Hence, we have $v'\in V(C)$ and $u'\in V(C')$.
This implies that $C'$ is a component of $\Gamma(\Omega',q)$.

If $u,u'\in V(C')$ are adjacent in $C'$ and $\mathcal{V}(u,C)=\{v_1,v_2,\ldots,v_k\}$, $k=|\mathcal{V}(u,C)|$, then there must be $k$ distinct vertices $v_1',v_2',\ldots,v_k'$ in $V(C)$
with $\Pi_{\Omega/\Omega'}(v_i')=u'$
such that $v_i'$ and $v_i$ are adjacent in $\Gamma(\Omega,q)$ for each $i$. Clearly, we have $\{v_1',v_2',\ldots,v_k'\}\subseteq \mathcal{V}(u',C)$
which implies $|\mathcal{V}(u',C)|\geq |\mathcal{V}(u,C)|$. Then, one can conclude that $t=|\mathcal{V}(u,C)|$ is independent of $u$.
Clearly, we see $1\leq t\leq q^{|\Omega|-|\Omega'|}$ and that $\Pi_{\Omega/\Omega'}$ is a $t$-to-1 graph homomorphism from $C$ to $C'$.\eproof

\medskip

We say sequence $\alpha\in\Omega$ is maximal in $\Omega$ if $\{\alpha0,\alpha1\}\cap\Omega=\emptyset$.
Clearly, for any maximal sequence $\alpha\in \Omega$, the set $\Omega'=\Omega\setminus\{\alpha\}$ also satisfies the condition C1.
For any component $C$ of $\Gamma(\Omega,q)$ and $u\in V(\Pi_{\Omega/\Omega'}(C))$, let $\rho(u,C)$ denote the set
of the entries indexed by $\alpha$ of the vertices in $\mathcal{V}(u,C)$, i.e.
\begin{align*}
\rho(u,C)=\left\{
\begin{array}{ll}
\{l_{\alpha}:[l]\in \mathcal{V}(u,C)\},&\text{if }u\in L(\Omega'),\\
\{r_{\alpha}:\dd{r}\in \mathcal{V}(u,C)\},&\text{if }u\in R(\Omega').
\end{array}
\right.
\end{align*}
For any component $C'$ of $\Gamma(\Omega',q)$, let $\mathbb{L}(C')$ denote the set of
components $C$ of $\Gamma(\Omega,q)$ with $C'=\Pi_{\Omega/\Omega'}(C)$.
The following theorem is a refinement of Theorem~\ref{gg01} for the case $\Omega'=\Omega\setminus\{\alpha\}$.

\begin{theorem}\label{gg02}
Suppose that $\alpha\in\Omega$ is a maximal sequence and $C'$ is a component of the graph $\Gamma(\Omega',q)$, where
$\Omega'=\Omega\setminus\{\alpha\}$.
Then, $s=|\mathbb{L}(C')|$ divides $q$ and there exist maps
$f: V(C')\to \mathbb{F}_q$, $g: \mathbb{L}(C')\to \mathbb{F}_q$ and an additive subgroup $G$ of $\mathbb{F}_q$ of order
$t=q/s$ such that, for any $C\in\mathbb{L}(C')$ and $u'\in V(C')$,
\begin{align}\label{cc03}
    \rho(u',C)=\left\{
    \begin{array}{ll}
    f(u')+g(C)+G,&\text{if }u'\in L(\Omega'),\\
    f(u')-g(C)+G,&\text{if }u'\in R(\Omega'),
    \end{array}\right.
\end{align}
where $\{g(C):C\in \mathbb{L}(C')\}$ is a representive set of the cosets of $G$ in $\mathbb{F}_q$, namely, $\{g(C)+G\}_{C\in\mathbb{L}(C')}$ are distinct cosets of $G$ in $\mathbb{F}_q$ with $\cup_{C\in\mathbb{L}(C')}\left(g(C)+G\right)=\mathbb{F}_q$.
Furthermore, if the component $C'$ contains the all-zero vertices,
then $t=|G|=p^{mn}$ for some integer $m$ with $0\leq m\leq \log_{p^{n}}q$, where $p$ is the characteristic of $\mathbb{F}_q$ and
$n>0$ is the smallest integer such that $q-1$ divides $(p^{n}-1)\gcd(q-1,w(\alpha)+1,|\alpha|-w(\alpha)+1)$.
\end{theorem}
{\bf Proof}:
Clearly, for any $x\in\mathbb{F}_q$, the automorphism $\theta_{x,\alpha}$ given in Lemma~\ref{lem:0002} defines a bijective map in each
$\mathbb{L}(C')$ and
\begin{align}
\label{t300}
\rho(u,\theta_{x,\alpha}(C))=\left\{
\begin{array}{ll}
\rho(u,C)+x,&\text{if }u\in L(\Omega'),\\
\rho(u,C)-x,&\text{if }u\in R(\Omega').
\end{array}
\right.
\end{align}
Since for any vertex $u\in V(C')$ the set $\{\rho(u,C):C\in\mathbb{L}(C')\}$ forms a partition of $\mathbb{F}_q$,
there is a unique component, denoted by $C_u$, in $\mathbb{L}(C')$ such that $0\in \rho(u,C_u)$.
Hence, for any $y\in \rho(u,C_u)$, we have
\begin{gather}
\label{t301}
\rho(u,C_u)=\left\{
\begin{array}{ll}
\rho(u,C_u)+y,&\text{if }u\in L(\Omega'),\\
\rho(u,C_u)-y,&\text{if }u\in R(\Omega'),
\end{array}
\right.\\
\mathcal{V}(u,C_u)=\theta_{y,\alpha}(\mathcal{V}(u,C_u)),
\end{gather}
and then
\begin{align}\label{t303}
C_u=\theta_{y,\alpha}(C_u).
\end{align}
From (\ref{t301}), we see that $\rho(u,C_u)$ is indeed an additive subgroup of $\mathbb{F}_q$.
Furthermore, from (\ref{t303}) we see that
$\rho(u',C_u)=\rho(u',C_u)+y$ is valid for any $u'\in V(C')$ and $y\in \rho(u,C_u)$.
Hence, $\rho(u',C_u)$ is a union of some cosets of $\rho(u,C_u)$ in $\mathbb{F}_q$.
Since according to Theorem~\ref{gg01} the cardinality of $\rho(u',C_u)$ is independent of $u'$ over $V(C')$,
the set $\rho(u',C_u)$ is just a coset of $\rho(u,C_u)$ in $\mathbb{F}_q$.
Furthermore, from (\ref{t300}) we see easily that
\begin{align*}
    \{\theta_{x,\alpha}(C):x\in\mathbb{F}_q\}=\mathbb{L}(C'), \text{ for any }C\in\mathbb{L}(C'),
\end{align*}
and there are some maps
$f: V(C')\to \mathbb{F}_q$ and $g: \mathbb{L}(C')\to \mathbb{F}_q$ such that (\ref{cc03}) is valid for the additive subgroup $G=\rho(u,C_u)$.
Since the order $t$ of the subgroup $G$ of $\mathbb{F}_q$ is a factor of $q$, we see that $s=|\mathbb{L}(C')|=q/t$ divides $q$.
Clearly, $\{g(C)|C\in \mathbb{L}(C')\}$ must be a representive set of the cosets of $G$ in $\mathbb{F}_q$.

Now we assume further that
the component $C'$ contains $u_0$, one of the all-zero vertices of $\Gamma(\Omega',q)$.
Since for any $x,y\in\mathbb{F}_q^*$ the automorphism $\lambda_{x,y}$ fixes the vertex $u_0\in V(C')$ and the component $C_{u_0}\in\mathbb{L}(C')$, we see that
\begin{align}\label{cc11}
    \xi^k\rho(u_0,C_{u_0})=\rho(u_0,C_{u_0}),
\end{align}
where $\xi$ is a primitive element of $\mathbb{F}_q$ and $k=\gcd(q-1,w(\alpha)+1,|\alpha|-w(\alpha)+1)$.
Since $G=\rho(u_0,C_{u_0})$ is an additive subgroup of $\mathbb{F}_q$,
we see that (\ref{cc11}) is equivalent to $xG=G$ for any nonzero element $x$ in
\begin{align*}
    \Delta=\left\{\left.\sum_{0\leq i<(q-1)/k}a_i\xi^{ki}\right|a_i\in\mathbb{F}_p\right\},
\end{align*}
where $p$ is the characteristic of $\mathbb{F}_q$.
Clearly, the set $\Delta$ is the smallest field which contains $\xi^k$ and thus
$\Delta=\mathbb{F}_{p^n}$, where $n>0$ is the smallest integer such that $(q-1)|(p^{n}-1)k$.
Hence, $G$ can be seen as a subspace of $\mathbb{F}_q$ over $\mathbb{F}_{p^n}$ and thus there is an integer $m$ with $0\leq m\leq \log_{p^n}q$
such that $t=|G|=p^{mn}$.
\eproof

\begin{corollary}
If $\mathbb{S}(\Omega)\subset\Omega$ and $\mathbb{S}(\Omega')\subset\Omega'$,
 then the number $t$ given in Theorem~\ref{gg01} is independent of the component $C$ of $\Gamma(\Omega,q)$.
Furthermore,
if $|\Omega|-|\Omega'|=1$, then the subgroup $G$ given in Theorem~\ref{gg02} is independent of the component $C'$ of $\Gamma(\Omega',q)$.
\end{corollary}
{\bf Proof}: Suppose $\mathbb{S}(\Omega)\subset\Omega$ and $\mathbb{S}(\Omega')\subset\Omega'$.
For any components $C_1$, $C_2$ of $\Gamma(\Omega,q)$, from Theorem~\ref{thm:-002} and the proof of Lemma~\ref{lem:0002}
we can see easily that there is an automorphism $\theta$ of $\Gamma(\Omega,q)$ such that $\theta(C_1)=C_2$
and, for any vertex $v\in V(C_1)$,
\begin{align*}
    \theta(\mathcal{V}(\Pi_{\Omega/\Omega'}(v),C_1))=\mathcal{V}(\Pi_{\Omega/\Omega'}(\theta(v)),C_2),
\end{align*}
which implies that the integer $t$ given in Theorem~\ref{gg01} is independent of the component $C$ of $\Gamma(\Omega,q)$.

Now we assume further that $|\Omega|-|\Omega'|=1$ and $\alpha$ is the sequence in $\Omega\setminus\Omega'$. Clearly, $\alpha$ is maximal in $\Omega$.
Let $\beta$ be an arbitrary sequence in $\Omega$ and $x$ an arbitrary element in $\mathbb{F}_q$.
For any component $C$ of $\Gamma(\Omega,q)$ and vertices $v, v'\in V(C)$ with $\Pi_{\Omega/\Omega'}(v)=\Pi_{\Omega/\Omega'}(v')$,
from the proof of Lemma~\ref{lem:0002} one can check easily that
the automorphism $\theta_{x,\beta}$ satisfies $\Pi_{\Omega/\Omega'}(\theta_{x,\beta}(v))=\Pi_{\Omega/\Omega'}(\theta_{x,\beta}(v'))$ and
$(\theta_{x,\beta}(v))_{\alpha}-v_{\alpha}=(\theta_{x,\beta}(v'))_{\alpha}-v'_{\alpha}$.
Then, from $(\theta_{x,\beta}(v))_{\alpha}-(\theta_{x,\beta}(v'))_{\alpha}=v_{\alpha}-v'_{\alpha}$,
the subgroup $G$ given in Theorem~\ref{gg02} for the component $C'=\Pi_{\Omega/\Omega'}(C)$ is the same one for the component $\Pi_{\Omega/\Omega'}(\theta_{x,\beta}(C))$.
Since any other component of $\Gamma(\Omega,q)$ is the image of $C$ under a few automorphisms of form $\theta_{x,\beta}$,
one can conclude further that the subgroup $G$ is independent of the component $C'$ of $\Gamma(\Omega',q)$.
\eproof

\section{Some Paths of $\Gamma(\Omega,q)$}\label{sec003}

In this section, we consider to express the vertices on a path by
their colors if the beginning vertex is of form $[0]_x$ or $\dd{0}_x$. However, we will deal with only the paths of $\Gamma(\Omega,q)$
which start at vertices in $L(\Omega)$. Similar argument is also true for the other case by symmetry.

Let
$[l^{(1)}]\dd{r^{(2)}}[l^{(3)}]\dd{r^{(4)}}\cdots$ be a given path
in $\Gamma(\Omega,q)$, namely, the vertices satisfy
\begin{gather}
\label{3002}
[l^{(2j-1)}]\neq [l^{(2j+1)}],\dd{r^{(2j)}}\neq\dd{r^{(2j+2)}},\\
\label{3001}
\{([l^{(2j-1)}],\dd{r^{(2j)}}),([l^{(2j+1)}],\dd{r^{(2j)}})\}\subset
E(\Gamma(\Omega,q)),
\end{gather}
for all $j\geq 1$. For convenience, we also write
\begin{align}
\label{q00}
x_i&=\left\{
\begin{array}{ll}
l^{(i)}_*, &\text{if }i\text{ is odd},\\
r^{(i)}_*, &\text{if }i\text{ is even}.
\end{array}\right.
\end{align}
Then, for $j\geq 1$, $s\geq 0$ and
$\alpha\in\{0,1\}^*\cup\{*\}$, from (\ref{3001}) and (\ref{000}) to
(\ref{002}), we have
\begin{gather}
\label{801}
    l^{(2j+1)}_{\alpha10^s}-l^{(2j-1)}_{\alpha10^s}=(x_{2j+1}-x_{2j-1})x_{2j}^sr^{(2j)}_{\alpha},
    \text{ if }\alpha10^s\in\Omega,\\
\label{802}
    r^{(2j+2)}_{\alpha01^s}-r^{(2j)}_{\alpha01^s}=(x_{2j+2}-x_{2j})x_{2j+1}^sl^{(2j+1)}_{\alpha},
    \text{ if }\alpha01^s\in\Omega.
\end{gather}

For positive integer $i$ and sequence $\boldsymbol{s}=(s_1,s_2,\ldots)$ of
nonnegative integers, let $\sigma_i(\boldsymbol{s})$ denote the sequence obtained from
$\mu_{1,i}(\boldsymbol{s})$ by deleting its first bit, namely,
\begin{align*}
\sigma_1(\boldsymbol{s})&=\left\{
\begin{array}{ll}
1^{s_1-1},&\text{if }s_1>0,\\
*,&\text{if }s_1=0,
\end{array}\right.
\end{align*}
and, for $i\geq 2$,
\begin{gather}
\label{t00'}
\sigma_{i}(\boldsymbol{s})=
\left\{
\begin{array}{ll}
\sigma_{i-1}(\boldsymbol{s})10^{s_{i}},&\text{if $i$ is even},\\
\sigma_{i-1}(\boldsymbol{s})01^{s_{i}},&\text{if $i$ is odd}.
\end{array}
\right.
\end{gather}

Clearly, for any sequence $\alpha\in\{0,1\}^*$,
there exist some nonnegative integers $s_1,s_2,\ldots,s_{2n}$ such that
\begin{align}
\alpha=1^{s_1}0^{s_2+1}1^{s_3+1}0^{s_4+1}\cdots1^{s_{2n-1}+1}0^{s_{2n}}
=\sigma_{2n}(s_1,s_2,\ldots,s_{2n}),
\end{align}
and consequently
\begin{align}
\alpha=\left\{\begin{array}{ll}
\sigma_{2n-1}(s_1,s_2,\ldots,s_{2n-2},s_{2n-1}+1),&\text{if }s_{2n}=0,\\
\sigma_{2n+1}(s_1,s_2,\ldots,s_{2n-1},s_{2n}-1,0),&\text{if }s_{2n}>0.
\end{array}\right.
\end{align}

Hereafter, let $\boldsymbol{s}=(s_1,s_2,\ldots)$ be a given sequence of nonnegative integers.
\begin{lemma}
If $[l^{(1)}]=[0]_{x_1}$,
then for any positive integer $i$,
\begin{gather}
\label{a01}
l^{(2i-1)}_{\sigma_{2j}(\boldsymbol{s})}=0, \text{ if }j\geq i\text{ and }\sigma_{2j}(\boldsymbol{s})\in\Omega,\\
\label{a02}
r^{(2i)}_{\sigma_{2j+1}(\boldsymbol{s})}=0, \text{ if }j\geq i\text{ and }\sigma_{2j+1}(\boldsymbol{s})\in\Omega.
\end{gather}
\end{lemma}
{\bf Proof}:
From $[l^{(1)}]=[0]_{x_1}$, we see that (\ref{a01}) is valid for $i=1$.
Furthermore, for $j\geq 1$ from $([l^{(1)}],\dd{r^{(2)}})\in E(\Gamma(\Omega,q))$ we have,
\begin{align*}
&r^{(2)}_{\sigma_{2j+1}(\boldsymbol{s})}
=r^{(2)}_{\sigma_{2j}(\boldsymbol{s})01^{s_{2j+1}}}+l^{(1)}_{\sigma_{2j}(\boldsymbol{s})01^{s_{2j+1}}}
=x_1^{s_{2j+1}}r^{(2)}_{\sigma_{2j}(\boldsymbol{s})0}
\\
=&x_1^{s_{2j+1}}\left(r^{(2)}_{\sigma_{2j}(\boldsymbol{s})0}+l^{(1)}_{\sigma_{2j}(\boldsymbol{s})0}\right)
=x_1^{s_{2j+1}}x_2l^{(1)}_{\sigma_{2j}(\boldsymbol{s})}=0.
\end{align*}
Hence, (\ref{a02}) is valid for $i=1$.

Now we assume that (\ref{a01}) and (\ref{a02}) are valid for some positive integer $i$.
For $j\geq i$, from (\ref{801}), (\ref{t00'}), (\ref{a01}) and (\ref{a02}),
\begin{align}
l^{(2i+1)}_{\sigma_{2j+2}(\boldsymbol{s})}
=l^{(2i+1)}_{\sigma_{2j+2}(\boldsymbol{s})}-l^{(2i-1)}_{\sigma_{2j+2}(\boldsymbol{s})}
=(x_{2i+1}-x_{2i-1})x_{2i}^{s_{2j+2}}r^{(2i)}_{\sigma_{2j+1}(\boldsymbol{s})}
=0.\label{a04}
\end{align}
From (\ref{802}), (\ref{t00'}), (\ref{a02}) and (\ref{a04}),
\begin{align*}
r^{(2i+2)}_{\sigma_{2j+3}(\boldsymbol{s})}
=r^{(2i+2)}_{\sigma_{2j+3}(\boldsymbol{s})}-r^{(2i)}_{\sigma_{2j+3}(\boldsymbol{s})}
=(x_{2i+2}-x_{2i})x_{2i+1}^{s_{2j+3}}l^{(2i+1)}_{\sigma_{2j+2}(\boldsymbol{s})}
=0.
\end{align*}
Hence, (\ref{a01}) and (\ref{a02}) are valid for $i+1$.

By induction, the proof is completed.\eproof\medskip

For positive integers $a_1,a_2,\ldots$ with
$a_{2j-1}\leq a_{2j}<a_{2j+1}$, $j=1,2,\ldots,$
let
\begin{align}
\label{a20}
\Phi_{\boldsymbol{s}}(a_1)=(x_{2a_1}-x_{2a_1-2})x_{2a_1-1}^{s_1},
\end{align}
where $x_0$ is defined as 0 if any.
For $n\geq 2$, let
\begin{align}
\Phi_{\boldsymbol{s}}(a_1,a_2,\ldots,a_{n})
=\left\{
\begin{array}{ll}
\Phi_{\boldsymbol{s}}(a_1,a_2,\ldots,a_{n-1})
(x_{2a_{n}+1}-x_{2a_{n}-1})x_{2a_{n}}^{s_{n}},&\text{if $n$ is even},\\
\Phi_{\boldsymbol{s}}(a_1,a_2,\ldots,a_{n-1})
(x_{2a_{n}}-x_{2a_{n}-2})x_{2a_{n}-1}^{s_{n}},&\text{if $n$ is odd}.
\end{array}\right.
\end{align}

\begin{theorem}
If $[l^{(1)}]=[0]_{x_1}$, then, for positive integers $i$, $j$ with $1\leq j\leq i$,
\begin{align}
\label{a06}
&r^{(2i)}_{\sigma_{2j-1}(\boldsymbol{s})}
=\sum_{1\leq a_1\leq a_2<a_3\leq a_4<\cdots<a_{2j-1}\leq i}
\Phi_{\boldsymbol{s}}(a_1,a_2,\ldots,a_{2j-1}),
\text{ if }\sigma_{2j-1}(\boldsymbol{s})\in\Omega,\\
\label{a07}
&l^{(2i+1)}_{\sigma_{2j}(\boldsymbol{s})}
=\sum_{1\leq a_1\leq a_2<a_3\leq a_4<\cdots<a_{2j-1}\leq a_{2j}\leq i}
\Phi_{\boldsymbol{s}}(a_1,a_2,\ldots,a_{2j}),\text{ if }\sigma_{2j}(\boldsymbol{s})\in\Omega.
\end{align}
\end{theorem}
{\bf Proof}:
For $1\leq j\leq i$, from (\ref{801}), (\ref{t00'}) and(\ref{a01}), we have
\begin{align}
\label{a10}
l^{(2i+1)}_{\sigma_{2j}(\boldsymbol{s})}
=\sum_{a=j}^i\left(l^{(2a+1)}_{\sigma_{2j}(\boldsymbol{s})}-l^{(2a-1)}_{\sigma_{2j}(\boldsymbol{s})}\right)
=\sum_{a=j}^i(x_{2a+1}-x_{2a-1})x_{2a}^{s_{2j}}r^{(2a)}_{\sigma_{2j-1}(\boldsymbol{s})}.
\end{align}
For $2\leq j\leq i$, from (\ref{802}), (\ref{t00'}) and (\ref{a02}), we have
\begin{align}
\label{a11}
r^{(2i)}_{\sigma_{2j-1}(\boldsymbol{s})}
=\sum_{a=j}^i\left(r^{(2a)}_{\sigma_{2j-1}(\boldsymbol{s})}-r^{(2a-2)}_{\sigma_{2j-1}(\boldsymbol{s})}\right)
=\sum_{a=j}^i(x_{2a}-x_{2a-2})x_{2a-1}^{s_{2j-1}}l^{(2a-1)}_{\sigma_{2j-2}(\boldsymbol{s})}.
\end{align}
If $s_1>0$, from $[l^{(1)}]=[0]_{x_1}$,
one can show easily that
\begin{align*}
r^{(2)}_{\sigma_1(s_1)}=x_2x_1^{s_1}=(x_2-x_0)x_1^{s_1},
\end{align*}
and then, for $i>1$, from $\sigma_1(s_1)=1^{s_1-1}=*01^{s_1-1}$ and (\ref{802}),
\begin{align*}
r^{(2i)}_{\sigma_1(s_1)}
=&r^{(2)}_{\sigma_1(s_1)}+\sum_{a=2}^i\left(r^{(2a)}_{\sigma_1(s_1)}-r^{(2a-2)}_{\sigma_1(s_1)}\right)\\
=&(x_2-x_0)x_1^{s_1}+\sum_{a=2}^i(x_{2a}-x_{2a-2})x_{2a-1}^{s_{1}-1}l^{(2a-1)}_{*}\\
=&\sum_{a=1}^i(x_{2a}-x_{2a-2})x_{2a-1}^{s_{1}}.
\end{align*}
Hence,  from $\sigma_1(0)=*$ and $x_0=0$, we see that
\begin{align}
\label{a12}
r^{(2i)}_{\sigma_1(s_1)}=&\sum_{a=1}^i(x_{2a}-x_{2a-2})x_{2a-1}^{s_{1}}
\end{align}
holds for any integers $s_1\geq 0$ and $i>0$, namely, (\ref{a06}) is valid for $j=1$.

Then, by induction on $j$, one can show (\ref{a06}) and (\ref{a07}) easily from the definition of $\Phi_{\boldsymbol{s}}(\cdot)$ and (\ref{a10}) to (\ref{a12}).
\eproof\medskip

From (\ref{a06}) and (\ref{a07}), we can also deduce the following corollary easily.
\begin{corollary}
If $[l^{(1)}]=[0]_{x_1}$, for $i\geq 1$ we have
\begin{align}
\label{a26}
r^{(2i)}_{\sigma_{2i-1}(\boldsymbol{s})}&=x_2x_1^{s_1}\prod_{a=2}^{2i-1}(x_{a+1}-x_{a-1})x_a^{s_{a}},
\text{ if }\sigma_{2i-1}(\boldsymbol{s})\in\Omega,\\
\label{a27}
l^{(2i+1)}_{\sigma_{2i}(\boldsymbol{s})}&=x_2x_1^{s_1}\prod_{a=2}^{2i}(x_{a+1}-x_{a-1})x_a^{s_{a}},
\text{ if }\sigma_{2i}(\boldsymbol{s})\in\Omega.
\end{align}
Furthermore, if $[l^{(1)}]=[0]_{x_1}$ and $x_2x_1^{s_1}=0$, for $i\geq 2$ we have
\begin{align}
\label{a28}
r^{(2i)}_{\sigma_{2i-3}(\boldsymbol{s})}&=\prod_{a=3}^{2i-1}(x_{a+1}-x_{a-1})x_a^{s_{a-2}},
\text{ if }\sigma_{2i-3}(\boldsymbol{s})\in\Omega,\\
\label{a29}
l^{(2i+1)}_{\sigma_{2i-2}(\boldsymbol{s})}&=\prod_{a=3}^{2i}(x_{a+1}-x_{a-1})x_a^{s_{a-2}},
\text{ if }\sigma_{2i-2}(\boldsymbol{s})\in\Omega.
\end{align}
\end{corollary}

\section{Existence of Some Cycles in $\Gamma(\Omega,q)$}\label{sec004}

In this section, we will show some conditions for the existence of
some cycles which contain a vertex of form $[0]_x$ or $\dd{0}_x$.
Therefore, some lower bounds of the girth of $\Gamma(\Omega,q)$ are deduced
from these conditions.

Let $C=[l^{(1)}]\dd{r^{(2)}}\cdots[l^{(2k-1)}]\dd{r^{(2k)}}$ be an arbitrary cycle in $\Gamma(\Omega,q)$ of length $2k$ with $[l^{(1)}]=[0]_{x_1}$.
For $1\leq j\leq 2k$, we define $x_j$ as in (\ref{q00}). Then, from (\ref{3001}) we see that
(\ref{3002}) is equivalent to
\begin{gather}
\label{a30}
x_j\neq x_{j+2}, 1\leq j\leq 2k,
\end{gather}
where $x_{2k+1}=x_1$ and $x_{2k+2}=x_2$.
We note that the cycle $C$ can also be expressed as $[l^{(1)}]\dd{r^{(2k)}}[l^{(2k-1)}]\cdots[l^{(3)}]\dd{r^{(2)}}$.

For $a\in\{0,1\}$, let $M_a$ denote the set of the sequences $\beta\in U$ that are lead by $a$, namely,
\begin{align*}
    M_a=\left\{
    \begin{array}{ll}
    \{0,01,010,0101,\ldots\},&\text{if }a=0,\\
    \{1,10,101,1010,\ldots\},&\text{if }a=1.
    \end{array}\right.
\end{align*}

\begin{lemma}
\label{lem:200}
Let $\beta$ be a sequence in $\Omega\cap U$.
\begin{enumerate}
  \item If $\beta\in M_0$, then $\Gamma(\Omega,q)$ has no cycle of length $2(|\beta|+1)$ containing a vertex of form $[0]_x$.
  \item If $\beta\in M_1$, then $\Gamma(\Omega,q)$ has no cycle of length $2(|\beta|+1)$ containing a vertex of form $\dd{0}_x$.
\end{enumerate}
\end{lemma}
{\bf Proof}:
We only prove the first conclusion. The second conclusion is valid by symmetry.

At first, we assume that $|\beta|$ is odd.
Clearly, $\beta=(01)^{i-1}0$ for some positive integer $i$. Assume that the cycle $C$ is of length $2k=4i=2(|\beta|+1)$.
For $0\leq s\leq 1$, from
$(01)^{i-1}0^s=\sigma_{2i}(0,\ldots,0,s)$ and (\ref{a27}) we have
\begin{align*}
    l^{(2i+1)}_{(01)^{i-1}0^s}&=x_{2i}^sx_2(x_3-x_1)\prod_{a=3}^{2i}(x_{a+1}-x_{a-1})\\
                       &=x_{2i+2}^sx_{4i}(x_{4i-1}-x_1)\prod_{a=3}^{2i}(x_{4i-a+1}-x_{4i-a+3}).
\end{align*}
If $l^{(2i+1)}_{(01)^{i-1}}=0$, then we have $x_2=0=x_{4i}$ which contradicts (\ref{a30}).
If $l^{(2i+1)}_{(01)^{i-1}}\neq0$, then we have
\begin{align*}
x_{2i}=l^{(2i+1)}_{(01)^{i-1}0}\left/ l^{(2i+1)}_{(01)^{i-1}}=x_{2i+2}\right.
\end{align*}
which still contradicts (\ref{a30}).

A proof for the case that $|\beta|$ is even can be given similarly by using (\ref{a26}).
\eproof

\begin{theorem}
\label{thm:y01}
Assume that $\mathbb{S}(\Omega)\subset\Omega$.
If $\beta\in\Omega\cap U$, then the girth of $\Gamma(\Omega,q)$ is at least $2(|\beta|+2)$.
\end{theorem}
{\bf Proof}:
This theorem is a simple corollary of Theorem~\ref{thm:-002} and Lemma~\ref{lem:200}.
\eproof
\medskip

\noindent {\bf Example 2}:
For positive integer $n$, let
\begin{align*}
\mathcal{X}_{n}&=\left\{
\begin{array}{ll}
\left(U_{4k+1}\setminus\{(10)^{k-1}1\}\right)\cup\{(01)^{k}0\},&\text{if }n=2k-1 \text{ is odd},\\
\left(U_{4k+3}\setminus\{(10)^k\}\right)\cup\{(01)^{k+1}\}, &\text{if }n=2k\text{ is even}.
\end{array}
\right.
\end{align*}
Clearly,
$\mathbb{S}(\Omega)\subset\Omega$ is valid for $\Omega=\mathcal{X}_n$,
$n\geq 1$. Hence, according to Theorem~\ref{thm:y01} we see easily that the girth of $\Gamma(\mathcal{X}_n,q)$ is at least $2n+8$.
We note that this lower bound of the girth of $\Gamma(\mathcal{X}_n,q)$ is
the same as the best known lower bound for the girth of $D(2n+3,q)=\Gamma(U_{2n+3},q)$
whose index set is as large as that of $\Gamma(\mathcal{X}_n,q)$.
\eproof

\begin{lemma}
\label{lem:210}
Suppose that $s,t$ are nonnegative integers with $\gcd(s-t,q-1)=1$.
\begin{enumerate}
  \item If $\beta\in M_0\cup\{\eta\}$ and $\{0^s\beta,0^t\beta\}\subset\Omega$, then $\Gamma(\Omega,q)$ has no cycle of length $2(|\beta|+2)$ containing a vertex of form $[0]_x$.
  \item If $\beta\in M_1\cup\{\eta\}$ and $\{1^s\beta,1^t\beta\}\subset\Omega$, then $\Gamma(\Omega,q)$ has no cycle of length $2(|\beta|+2)$ containing a vertex of form $\dd{0}_x$.
\end{enumerate}
\end{lemma}
{\bf Proof}: As in Lemma~\ref{lem:200}, we only prove the first conclusion of this lemma.

At first, we assume that $|\beta|$ is even.
Clearly, $\beta=(01)^{i-1}$ for some positive integer $i$. Assume that the cycle $C$ is of length $2k=4i=2(|\beta|+2)$.
From $0^{s}(01)^{i-1}=\sigma_{2i}(0,s,0,\ldots,0)$ and (\ref{a27}) we have
\begin{align*}
    l^{(2i+1)}_{0^s(01)^{i-1}}&=x_{2}^{s+1}(x_3-x_1)\prod_{a=3}^{2i}(x_{a+1}-x_{a-1})\\
    &=x_{4i}^{s+1}(x_{4i-1}-x_1)\prod_{a=3}^{2i}(x_{4i-a+1}-x_{4i-a+3}).
\end{align*}
Since these two equalities are still true when $s$ is replaced by $t$,
from $\gcd(s-t,q-1)=1$ we can get easily $x_2=x_{4i}$ which contradicts (\ref{a30}).

A proof for the case that $|\beta|$ is odd can be given similarly by using (\ref{a26}).
\eproof

\begin{theorem}
\label{thm:y02}
Assume that $\mathbb{S}(\Omega)\subset\Omega$.
If $a$ is a symbol in $\{0,1\}$ and $\beta$ is a sequence in $M_a\cup\{\eta\}$ such that $\{a^s\beta,a^t\beta\}\subset\Omega$ for some nonnegative integers $s,t$ with $\gcd(s-t,q-1)=1$, then the girth of $\Gamma(\Omega,q)$ is at least $2(|\beta|+3)$.
\end{theorem}
{\bf Proof}:
This theorem is a simple corollary of Theorem~\ref{thm:-002} and Lemma~\ref{lem:210}.
\eproof
\medskip

\noindent {\bf Example 3}:
Let
\begin{align*}
    \Omega_3=\{\eta,0,1,01,10,010,0101,01010,0^2,0^21,0^210,0^2101,0^21010\}.
\end{align*}
From $\mathbb{S}(\Omega_3)\cup\{01010,0^21010\}\subset\Omega_3$
and Theorem~\ref{thm:y02}, the girth of $\Gamma(\Omega_3,q)$ is at least $2(5+3)=16$.
\eproof

\begin{lemma}
\label{lem:h00}
\begin{enumerate}
  \item If $\beta\in M_0\cup\{\eta\}$ and $1^m\beta\in\Omega$ for some positive integer $m$,
  then $\Gamma(\Omega,q)$ has no cycle of length $2(|\beta|+3)$ containing the vertex $[0]_0$.
  \item If $\beta\in M_1\cup\{\eta\}$ and $0^m\beta\in\Omega$ for some positive integer $m$,
  then $\Gamma(\Omega,q)$ has no cycle of length $2(|\beta|+3)$ containing the vertex $\dd{0}_0$.
\end{enumerate}
\end{lemma}
{\bf Proof}: As in Lemma~\ref{lem:200}, we only prove the first conclusion of this lemma.

At first, we assume that $\beta$ is of even length,
namely, $\beta=(01)^{i-1}$ for some positive integer $i$. Assume that the cycle $C$ is of length $2k=4i+2=2(|\beta|+3)$ and the color $x_1$ of the vertex $[l^{(1)}]=[0]_{x_1}$ is 0.

For $0\leq s\leq 1$,
from $\sigma_{2i-1}(m,0,\ldots,0,s)=1^{m}(01)^{i-2}01^s\in\Omega$ and (\ref{a28}), we have
\begin{align*}
r^{(2i+2)}_{\sigma_{2i-1}(m,0,\ldots,0,s)}
&=x_3^{m}x_{2i+1}^s\prod_{a=3}^{2i+1}(x_{a+1}-x_{a-1})\\
&=x_{4i+1}^{m}x_{2i+3}^s\prod_{a=3}^{2i+1}(x_{4i-a+3}-x_{4i-a+5}),
\end{align*}
where $\sigma_{2i-1}(m,0,\ldots,0,s)$ denotes the sequence $\sigma_{1}(m+s)=1^{m-1+s}$ when $i=1$.

From $x_1=0$ we see $x_3\neq 0$ and $r^{(2i+2)}_{\sigma_{2i-1}(m,0,\ldots,0)}\neq 0$, and thus we have
\begin{align*}
x_{2i+1}=r^{(2i+2)}_{\sigma_{2i-1}(m,0,\ldots,0,s)}\left/r^{(2i+2)}_{\sigma_{2i-1}(m,0,\ldots,0)}=x_{2i+3}\right.
\end{align*}
which contradicts (\ref{a30}).

A proof for the case that $\beta$ is of odd length can be given similarly by using (\ref{a29}).
\eproof

\begin{theorem}
\label{thm:800}
Let $a$ be a symbol in $\{0,1\}$ and $\bar{a}$ the other.
Assume that $\mathbb{S}(\Omega)\cup\mathbb{H}_{a}(\Omega)\subset\Omega$.
 If there are positive integers $t, n_1, n_2, \ldots, n_t, m$ and sequences $\gamma,\beta\in M_{a}\cup\{\eta\}$ with $|\beta|=|\gamma|+2t-1$ such that
 $\{a\bar{a}^{n_1}a\bar{a}^{n_2}\cdots a\bar{a}^{n_t}\gamma,\bar{a}^m\beta\}\subset\Omega$,
 then the girth of $\Gamma(\Omega,q)$ is at least $2(|\beta|+4)$.
\end{theorem}
{\bf Proof}:
Without loss of generality we assume that $a=1$.
Suppose $\mathbb{S}(\Omega)\cup\mathbb{H}_{1}(\Omega)\subset\Omega$ and
$\{10^{n_1}10^{n_2}\cdots10^{n_t}\gamma,0^m\beta\}\subset\Omega$
 for some positive integers $t, n_1, n_2, \ldots, n_t, m$ and sequences $\gamma, \beta\in M_{1}\cup\{\eta\}$ with $|\beta|=|\gamma|+2t-1$.
Assume that the cycle $C$ is of length $2k=2(|\beta|+3)$.
From Theorem~\ref{thm:-002} and Lemma~\ref{lem:h00}, we can deduce easily $x_{2i}\neq 0$ for $i=1,2,\ldots,k$.
According to Theorem~\ref{thm:401}, we can assume without loss of generality that the color $x_1$ of $[l^{(1)}]=[0]_{x_1}$ is 0.

Now we assume that $\beta$ is of even length. Then, there is a positive integer $j$ such that
$\gamma=(10)^{j-1}1$ and $k=2t+2j+1$.
Let $\omega=10^{n_1}10^{n_2}\cdots10^{n_t}(10)^{j-1}$.
For $0\leq s\leq 1$, from $\omega1^s=\sigma_{2t+2j-1}(1,n_1-1,0,n_2-1,\ldots,0,n_t-1,0,\ldots,0,s)\in\Omega$ and (\ref{a28}) we see
\begin{align*}
    r^{(2t+2j+2)}_{\omega1^s}&=x_3x_{2t+2j+1}^s\prod_{i=1}^tx_{2i+2}^{n_i-1}\prod_{a=3}^{2t+2j+1}(x_{a+1}-x_{a-1})\\
    &=x_{4t+4j+1}x_{2t+2j+3}^s\prod_{i=1}^tx_{4t+4j+2-2i}^{n_i-1}\prod_{a=3}^{2t+2j+1}(x_{4t+4j-a+3}-x_{4t+4j-a+5}).
\end{align*}
Since from $x_1=0$ we have $x_3\neq 0$, then we see easily $r^{(2t+2j+2)}_{\omega}\neq 0$ and
\begin{align*}
    x_{2t+2j+1}=r^{(2t+2j+2)}_{\omega1}\left/r^{(2t+2j+2)}_{\omega}=x_{2t+2j+3}\right.
\end{align*}
which contradicts (\ref{a30}). Hence, $\Gamma(\Omega,q)$ has no cycle of length $2k=2(|\beta|+3)$ if $\beta$ is of even length.

By using (\ref{a29}), one can show similarly that $\Gamma(\Omega,q)$ has no cycle of length $2k=2(|\beta|+3)$ if $\beta$ is of odd length.

Furthermore, if we replace the sequences
$10^{n_1}10^{n_2}\cdots10^{n_t}\gamma,0^m\beta$ in the above argument by their subsequences obtained by deleting a few rightmost bits,
it can be concluded that $\Gamma(\Omega,q)$ has no cycle of length between $8$ and $2(|\beta|+3)$.
On the other hand, from $10\in \Omega$ and Theorem~\ref{thm:y01}, $\Gamma(\Omega,q)$ has no cycle of length less than 8.
Hence, the girth of $\Gamma(\Omega,q)$ is at least $2(|\beta|+4)$.\eproof
\medskip

\noindent{\bf Example 4:}
Let
\begin{align*}
    \Omega_4=\{10^310,   10^31,10^3,10^2,10,1;
    0^2101,   0^210,0^21,0^2,0;
    01;
    0^310,   0^31,0^3;
    \eta\}.
\end{align*}
From $\mathbb{S}(\Omega_4)\cup\mathbb{H}_{1}(\Omega_4)\cup\{10^310,0^2101\}\subset\Omega_4$
and Theorem~\ref{thm:800}, the girth of $\Gamma(\Omega_4,q)$ is at least $2(3+4)=14$.
\eproof

\begin{theorem}
\label{cor009} Assume that $\mathbb{H}_0(\Omega)\cup\mathbb{H}_1(\Omega)\subset\Omega$.
If $\alpha\in\Omega\cap U$, then the girth of $\Gamma(\Omega,q)$ is at least $2(|\alpha|+3)$.
\end{theorem}
{\bf Proof}:
Since $\mathbb{H}_0(\Omega)\cup\mathbb{H}_1(\Omega)\subset\Omega$,
according to Theorem~\ref{thm:-001} we can only deal with the cycles containing the edge $([0]_0, \dd{0}_0)$.
Clearly, $\Omega$ must contain all the sequences in $U$ which
are shorter than $\alpha$.
Then, according to Lemma~\ref{lem:h00}, we see that $\Gamma(\Omega,q)$ has no cycle of length less than $2(|\alpha|+3)$.
Hence, the girth of $\Gamma(\Omega,q)$ is at least $2(|\alpha|+3)$.
\eproof\medskip

Therefore we can deduce easily the following corollary, which was first proved in \cite{L.U.1995}.
\begin{corollary}
\label{cor:201}
For $k\geq 2$, the girth of $D(k,q)$ is at least $k+4$.
\end{corollary}
{\bf Proof}:
Since $\mathbb{H}_0(\Omega)\cup\mathbb{H}_1(\Omega)\subset\Omega$ is valid for $\Omega=U_k$ and the maximum length of the sequences in $U_k$
is $\lfloor (k-1)/2\rfloor$, according to
Theorem~\ref{cor009} the girth of $D(k,q)=\Gamma(U_k,q)$ is not smaller than
$2\lfloor (k-1)/2\rfloor+6\geq k+4$.
\eproof
\medskip

We note that Theorem~\ref{cor009} can also be deduced simply from Corollary~\ref{cor:201}.
Indeed, if $\mathbb{H}_0(\Omega)\cup\mathbb{H}_1(\Omega)\subset\Omega$ and $\alpha\in\Omega\cap U$,
without loss of generality we assume further $\alpha\in M_0$, then we have $U_k\subset \Omega$ for $k=2|\alpha|+1$
and thus the girth of $\Gamma(\Omega,q)$ is not smaller than that of $\Gamma(U_k,q)=D(k,q)$.
Hence, from Corollary~\ref{cor:201} and $k+4=2|\alpha|+5$ we see that the girth of $\Gamma(\Omega,q)$ is at least $2(|\alpha|+3)$.

\section{Conclusion}\label{sec005}
To generalize the bipartite graph
$D(k,q)$ proposed by Lazebnik and Ustimenko,
we construct in this paper a bipartite graph $\Gamma(\Omega,q)$ for any set $\Omega$ of binary sequences
that are employed to index the
entries of the vertex vectors.
Sufficient conditions for the generalized graph $\Gamma(\Omega,q)$ to admit a variety of automorphisms are proposed.
A sufficient condition for $\Gamma(\Omega,q)$ to be edge-transitive is shown by using these automorphisms.
For $\Gamma(\Omega,q)$, we show some invariants which show that $\Gamma(\Omega,q)$ is disconnected in general.
For the paths and cycles which contain a vertex of form $[0]_x$ or $\dd{0}_x$, we show an expression for each vertex on them in terms of the colors of the vertices. From these expressions, we deduce a few lower bounds
for the girth of $\Gamma(\Omega,q)$. We note that the results obtained in this paper generalize many of the known results on $D(k,q)$.
Furthermore, one can propose easily some conditions for the generalized graphs to be a family
of graphs with large girth in the sense proposed by Biggs. For example, the graphs $\Gamma(\mathcal{X}_n,q)$ form a such family.



\end{document}